\documentclass[preprint,11pt]{elsarticle}
\usepackage[bookmarks=true,colorlinks=true,linkcolor=blue]{hyperref}

\pdfoutput=1

\biboptions{sort&compress}

\usepackage{color}
\usepackage{amssymb,amsmath}
\usepackage{amsthm}

\newtheorem{theorem}{Theorem}
\newtheorem{algorithm}{Algorithm}

\usepackage{calc}

\usepackage{calc}
\usepackage[margin=1.in]{geometry}

\renewcommand{\url}[1]{}
\newcommand{\citeCount}[1]{}

\newtheorem{condition}{Condition}
\newtheorem{observation}{Observation}

\newtheorem*{ModelProblemI}{Model Problem MP-I1}
\newtheorem*{ModelProblemVa}{Model Problem MP-V1}
\newtheorem*{ModelProblemVb}{Model Problem MP-V2}
\newcommand{\dx}{\Delta x}
\newcommand{\dy}{\Delta y}
\newcommand{\dt}{\Delta t}

\newcommand{\bogus}[1]{{}}

\newcommand{\ev}{\mathbf{ e}}

\newcommand{\iv}{\mathbf{ i}}

\newcommand{\nv}{\mathbf{ n}}

\newcommand{\qv}{\mathbf{ q}}

\newcommand{\uv}{\mathbf{ u}}
\newcommand{\vv}{\mathbf{ v}}

\newcommand{\xv}{\mathbf{ x}}

\newcommand{\Fv}{\mathbf{ F}}

\newcommand{\Iv}{\mathbf{ I}}

\newcommand{\Lv}{\mathbf{ L}}

\newcommand{\half}{{1\over2}}

\newcommand{\Real}{{\mathbb R}}
\newcommand{\Complex}{{\mathbb C}}

\newcommand{\tableFont}{\scriptsize}

\newcommand{\Bc}{{\mathcal B}}

\newcommand{\Lc}{{\mathcal L}}

\newcommand{\tauv}{\boldsymbol{\tau}}

\newcommand{\sigmav}{\boldsymbol{\sigma}}

\newcommand{\grad}{\nabla}

\newcommand{\tableFontSize}{\footnotesize}
\newcommand{\num}[2]{#1e#2} 

\newcommand{\rateLabel}{rate}

\newcommand{\rhos}{\bar{\rho}}

\newcommand{\us}{\bar{u}}

\newcommand{\usv}{\bar{\uv}}
\newcommand{\vsv}{\bar{\vv}}

\newcommand{\vvs}{\bar{\vv}}
\newcommand{\uvs}{\bar{\uv}}

\newcommand{\amp}{A}

\newcommand{\normalss}{\sffamily}

\newcommand{\kx}{k_x}

\newcommand{\Lt}{\Lc}

\newcommand{\OmegaF}{\Omega^F}
\newcommand{\OmegaFh}{\Omega^F_h}
\newcommand{\OmegaS}{\Omega^S}

\newcommand{\etah}{\hat{\eta}}
\newcommand{\usvh}{\hat{\usv}}
\newcommand{\uvh}{\hat{\uv}}
\newcommand{\vvh}{\hat{\vv}}
\newcommand{\vh}{\hat{v}}
\newcommand{\ph}{\hat{p}}
\newcommand{\Lsv}{\bar{\Lv}}
\newcommand{\Ks}{\bar{K}}
\newcommand{\Ts}{\bar{T}}
\newcommand{\Bs}{\bar{B}}
\newcommand{\hs}{\bar{h}}

\newcommand{\GammaI}{\Gamma}
\newcommand{\GammaIh}{\Gamma_h}
\newcommand{\Hf}{H}
\newcommand{\hf}{h_f}

\newcommand{\vtz}{\tilde{v}}
\newcommand{\ustz}{\tilde{u}}
\newcommand{\ptz}{\tilde{p}}
\newcommand{\fx}{f_x}
\newcommand{\ft}{f_t}
\newcommand{\dsf}{\delta}

\newcommand{\Ck}{C_k}
\newcommand{\Sk}{S_k}
\newcommand{\Ca}{C_\alpha}
\newcommand{\Sa}{S_\alpha}

\newcommand{\tanv}{\ev}

\newcommand{\expkw}{e^{i(k x -\omega t)}}
\newcommand{\ampe}{\us_{\rm max}}

\newcommand{\strutt}{\rule{0pt}{10pt}}

\newcommand{\len}{{\bar\ell}}

\usepackage{tikz}

\newlength{\tfwidth}
\newlength{\tfheight}
\newlength{\tfxa}
\newlength{\tfxb}
\newlength{\tfya}
\newlength{\tfyb}
\newcommand{\trimFigNoBox}[6]{%
\setlength\fboxsep{1pt}
\setlength\fboxrule{0.0pt}
\fbox{\includegraphics[width=#2, clip, trim=#3 #4 #5 #6]{#1}}%
}

\newcommand{\trimFig}[6]{%
\setlength{\tfwidth}{(#2+#2*\real{#3})+#2*\real{#4}}
\setlength{\tfheight}{(#2+#2*\real{#5})+#2*\real{#6}}%
\setlength{\tfxa}{\tfwidth*\real{#3}}%
\setlength{\tfxb}{\tfwidth*\real{#4}}%
\setlength{\tfya}{\tfheight*\real{#5}}%
\setlength{\tfyb}{\tfheight*\real{#6}}%
\trimFigNoBox{#1}{#2}{\tfxa}{\tfya}{\tfxb}{\tfyb}%
}
\begin{document}

\small

\begin{frontmatter}
\title{A stable partitioned FSI algorithm for incompressible flow\\ and structural shells}

\author[llnl]{J. W. Banks\fnref{llnlThanks}}
\ead{banks20@llnl.gov}

\author[llnl]{W. D. Henshaw\corref{cor1}\fnref{llnlThanks}}
\ead{henshaw1@llnl.gov}

\author[rpi]{D. W. Schwendeman\fnref{donThanks}}
\ead{schwed@rpi.edu}

\address[llnl]{Centre for Applied Scientific Computing, Lawrence Livermore National Laboratory, Livermore, CA 94551, USA}

\address[rpi]{Department of Mathematical Sciences, Rensselaer Polytechnic Institute, Troy, NY 12180, USA}

\cortext[cor1]{Corresponding author. Mailing address: Centre for Applied Scientific Computing, L-422, Lawrence Livermore National Laboratory, Livermore, CA 94551, USA. Phone: 925-423-2697. Fax: 925-424-2477. }

\fntext[llnlThanks]{This work was performed under the auspices of the U.S. Department of Energy (DOE) by
  Lawrence Livermore National Laboratory under Contract DE-AC52-07NA27344 and by 
  DOE contracts from the ASCR Applied Math Program.}

\fntext[donThanks]{This research was supported by Lawrence Livermore National Laboratory under
subcontract B548468, and by the National Science Foundation under grant DMS-1016188.}

\begin{abstract}
A stable partitioned algorithm for fluid-structure
interaction (FSI) problems that couple viscous incompressible flow with structural
shells or beams is described.
This added-mass partitioned (AMP) scheme uses Robin (mixed) interface
conditions for the pressure and velocity in the fluid that are derived
directly from the governing equations.  The AMP scheme is stable even for very light structures, 
requires no sub-iterations, and can be made fully second-order, or higher-order, accurate.  
The new scheme is shown to be stable through the analysis of a model problem. 
Exact traveling wave solutions for three FSI model problems are derived. 
Numerical results for a linearized FSI problem in two-dimensions, using these exact solutions,
demonstrate that the
scheme is stable and accurate, even for very light structures.  

\end{abstract}

\begin{keyword}
fluid-structure interaction, added mass instability, incompressible fluid flow, structures, shells, beams
\end{keyword}

\end{frontmatter}

\tableofcontents

\section{Introduction}

We consider fluid-structure interaction (FSI) problems that couple an 
incompressible fluid with a structural shell or beam.
Partitioned schemes (also known as modular or sequentitial schemes)
solve the FSI problem numerically by 
splitting the solution of the equations into separate solvers for the fluid and structual domains.
These are in contrast to monolithic schemes whereby the fluid and solid equations are solved together
as a large system of equations. 
Strongly coupled partitioned schemes perform multiple sub-iterations per time-step to solve the
coupled equations while loosely coupled schemes use only a few or no iterations.
The traditional partitioned algorithm for shells uses the velocity of the solid as a boundary condition on the fluid.
The force of the fluid is accounted for through a body forcing on the shell. 
It has been found that partitioned schemes may be unstable, or require multiple sub-iterations per time-step, 
when the density of the structure is similar to or lighter than that of the fluid~\cite{CausinGerbeauNobile2005,vanBrummelen2009}.
These instabilities are attributed to the {\em added-mass effect} whereby the force required to accelerate a structure immersed
in a fluid must also account for accelerating the surrounding fluid.
The added-mass effect has been found to be
 especially problematic in many biological flows such as haemodynamics since the
density of the fluid (blood) is similar to that of the adjacent structure (arterial walls)~\cite{YuBaekKarniadakis2013}. 

In this work we develop a new stable partitioned algorithm for coupling
incompressible flows with structural shells (or beams) that overcomes the added-mass
effects.
The scheme is developed and evaluated for a linearized problem where the fluid
is modeled with the Stokes equations on a fixed reference domain
and the structure is modeled with a linear beam or generalized string model.
This added-mass partitioned (AMP) scheme is stable and requires no sub-iterations.
We
do, however, generally use a predictor-corrector algorithm for the fluid with
one correction step, since this increases the stability region of the
incompressible flow solver.  The AMP approximation can be made fully
second-order (or high-order) accurate.  The approach uses generalized Robin-type
(or ``mixed'') boundary conditions at the interface where the coupling
coefficients are derived directly from the governing equations.
The key coupling equation is a
Robin condition on the fluid pressure that is formed by combining the evolution
equation for the shell with the fluid momentum equation on the interface. This
procedure also provides a generalized Robin condition for the tangential
components of the velocity.  For heavy solids the AMP scheme reduces to the
traditional coupling (velocity defined from the structure), while for heavy fluids the scheme
approaches a free surface problem (traction defined from the structure).
The stability of the new AMP scheme is proved for a two-dimensional
model problem.  The incompressible flow equations in our FSI scheme are solved
with a fractional-step (split-step) method based on the velocity-pressure
formulation~\cite{ICNS,splitStep2003} that can be made fully second-order (or
higher-order) accurate.
We develop exact traveling
wave solutions to three FSI model problems and these are used to verify that the
AMP scheme is stable and second-order accurate in the maximum-norm.

Previously, a stable partitioned algorithm for compressible flows and elastic (bulk) solids that overcomes the
added-mass effect was developed in Banks et~al.~\cite{fsi2012} and used deforming composite grids
to treat large solid motions. 
For compressible flows the added-mass effect is more localized than for incompressible flows~\cite{vanBrummelen2009}
due to the finite speeds of wave propagation in the fluid. 
The scheme in~\cite{fsi2012} was based on the analysis
of Banks and Sj\"ogreen~\cite{sjogreenBanks2012}. The approach uses a local characteristic analysis
of a fluid-structure Riemann problem to define  an impedance weighted averaging of the fluid and
solid interface values (i.e.~a Robin-Robin coupling). 
The algorithm was extended to the case of rigid bodies
in~\cite{lrb2013} where it was shown that the scheme remains stable even for rigid bodies of zero mass.


In recent work~\cite{fib2013r}, we have also developed a stable partitioned FSI algorithm for
incompressible flow coupled to bulk elastic solids that remains
stable even for light structures when added mass effects are large. The case of bulk solids
is different from the situation of shells or beams since the fluid traction 
enters as a boundary condition on the bulk solid rather than as a body forcing.
The AMP algorithm for bulk solids is thus different than the one for shells or beams
presented in this article, although both algorithms impose generalized Robin boundary conditions
on the fluid domain.



Problems involving FSI is a very active field of
research with many publications, see for example~\cite{DoneaGiulianiHalleux1982,CebralLohner1997,Lohner99,LeTallecMouro2001,PipernoFarhat2001,SchaferTeschauer2001,GlowinskiPanHelsaJosephPeriaux2001,Guruswamy2002,MillerColella2002,KuhlHulshoffBorst2003,ArientiHungMoranoShepherd2003,Farhat2006,HronTurekMonolithic2006,Tezduyar2006,SchaferHeckYigit2006,AhnKallinderis2006,CirakDeiterdingMauch2007,vanLoon2007,BorazjaniGeSotiropolous2008,Wilcox2010,BartonObadiaDrikakis2011,DegrooteVierendeels2011,TezduyarTakizawaBrummerChen2011,CrosettoDeparisFouresteyQuarteroni2011,Degroote2011,GretarssonKwatraFedkiw2011,HouWangLayton2012}. There has also
been much work concerning partitioned algorithms and the added-mass effect and we outline some of this work now.  
%
For the case of incompressible fluids coupled to thin structural shells, as considered in this article,
there has been some success in 
treating the added-mass instability. 
The kinematically coupled 
scheme of Guidoboni et~al.~\cite{GuidoboniGlowinskiCavalliniCanic2009}, later extended to the
$\beta$-scheme by \v{C}ani\'c, Muha and Buka\v{c}~\cite{CanicMuhaBukac2012}, is a stable partitioned scheme
that uses operator splitting. Nobile and Vergara~\cite{NobileVergara2008}, among others, have developed 
stable semi-monolithic schemes for the case of thin shells. 
Even fully monolithic schemes for thin shells~\cite{FigueroaVignonClementelJansenHughesTaylor2006} are
not overly expensive when the number of degrees of freedom for the structural shell is small compared
to those in the fluid. However, there are still numerous advantanges to developing a modular scheme, such as the one discussed here, since, for example, existing fluid
and solid solvers can be used without major changes. 

Numerous authors have analyzed the added-mass effect and the stability of FSI
algorithms~\cite{Conca1997,LeTallecMani2000,MokWallRamm2001,ForsterWallRamm2007,AstorinoChoulyFernandez2009,vanBrummelen2009,Gerardo_GiordaNobileVergara2010}.
Causin, Gerbeau and Nobile~\cite{CausinGerbeauNobile2005}, for example, analyze a model
problem of a structural shell coupled to an incompressible fluid and show that
the traditional scheme can be unconditionally unstable over a certain range of
parameters.
Many partitioned algorithms require multiple sub-iterations per time-step to overcome the added-mass instability, 
and there have been a number of approaches that have been developed to reduce the number of sub-iterations~\cite{MokWallRamm2001}. 
Robin interface conditions are used to stabilize partitioned schemes for both shells and bulk solids, by, for example,
Nobile, Vergara and co-workers~\cite{NobileVergara2008,BadiaNobileVergara2008,BadiaNobileVergara2009,Gerardo_GiordaNobileVergara2010,NobileVergara2012}, and Astorino, Chouly and Fernandez~\cite{AstorinoChoulyFernandez2009}.
Baek and Karniadakis~\cite{BaekKarniadakis2012} and Yu, Baek and Karniadakis~\cite{YuBaekKarniadakis2013} have developed fictitious pressure and
fictitious mass algorithms which incorporate additional terms into the governing equations to account
for added-mass effects. 
%
Degroote et~al.~\cite{DegrooteSwillensBruggemanHaeltermanSegersVierendeels2010} have developed an interface artificial compressibility method
that adds a source term to the fluid continuity equation near the interface; the source term goes away as the sub-iterations converge.
Idelsohn et~al.~\cite{IdelsohnDelPinRossiOnate2012}, and Badia, Quaini and Quarteroni~\cite{BadiaQuainiQuarteroni2008} 
form approximate factorizations of the fully
monolithic scheme to construct partitioned schemes but these still may require many iterations to converge. 
Degroote et~al.~\cite{DegrooteBruggemanHaeltermanVierendeels2008,DegrooteBatheVierendeels2009}
use reduced order models and Aitken acceleration methods to reduce the number of
iterations in partitioned schemes.

The remainder of the manuscript is organized as follows. In Section~\ref{sec:governingEquations} we define the governing equations.
The AMP interface conditions are derived in Section~\ref{sec:AMPalgorithm}, and a second-order accurate predictor-corrector
algorithm based on these conditions is described.
In Section~\ref{eq:analysis} we analyze the stability of the AMP scheme for an inviscid  model problem.
Some FSI model problems are defined in~\ref{sec:modelProblems}.  
Numerical results for the model problems are presented in Section~\ref{sec:results}, using
both the method of analytic solutions as well as solving for some exact traveling wave solutions.
Conclusions are provided in Section~\ref{sec:conclusions}.
The exact traveling wave solutions
to the FSI model problems are presented in~\ref{sec:travelingWave}.

\section{Governing equations for incompressible flow and a structural shell} \label{sec:governingEquations}

Consider a fluid-structure interaction problem in which
an incompressible viscous fluid in a two-dimensional domain $\OmegaF$ is coupled to a structural shell (or beam) on the interface $\GammaI$, 
which is a smooth curve on a portion of the boundary of $\OmegaF$.
Since, for the purposes of this article, the primary concern is the stability of numerical algorithms, 
we consider the situation of small perturbations to an equilibrium state. 
Thus the nonlinear advection terms in the Navier-Stokes equations are
neglected and the fluid domain $\OmegaF$ (including the interface $\GammaI$) is fixed in time. The equations governing the structural shell
are taken to be a linearized beam or generalized string model. 
The fluid, therefore, is governed by the Stokes equations on a fixed reference domain $\OmegaF$,
which in velocity-divergence form are given by,
\begin{alignat}{3}
  &  \rho \frac{\partial\vv}{\partial t} 
                 =  \grad\cdot\sigmav  , \qquad&& \xv\in\OmegaF ,  \label{eq:NS3dv} \\
  & \grad\cdot\vv =0,  \quad&& \xv\in\OmegaF ,  \label{eq:NS3dDiv}
\end{alignat}
where $\rho>0$ is the constant density and $\vv=\vv(\xv,t)$ is the velocity vector at a position $\xv$ and time $t$. 
The fluid stress tensor, $\sigmav$, is given by 
\begin{align*}
&  \sigmav = -p \Iv + \tauv, \qquad
  \tauv  = \mu [ \grad\vv + (\grad\vv)^T ], 
\end{align*}
where $p=p(\xv,t)$ is the pressure, $\Iv$ the identity tensor, $\mu\ge 0$ the constant fluid viscosity and $\tauv$ the viscous stress tensor.
These equations require initial conditions for $\vv(\xv,0)$ and appropriate boundary conditions.
For future reference, 
the components of a vector, such as $\vv$ will be denoted by $v_m$, $m=1,2,3$, (i.e. $\vv=[v_1, v_2, v_3]^T$), while components
of a tensor such as $\sigmav$, will be denoted by $\sigma_{mn}$, $m,n=1,2,3$.
The fluid equations can also be written in velocity-pressure form, 
\begin{alignat}{3}
  & \rho \frac{\partial\vv}{\partial t} 
                 =  \grad\cdot\sigmav  , \qquad&& \xv\in\OmegaF ,  \label{eq:VPv} \\
  & \Delta p = 0,  \quad&& \xv\in\OmegaF , \label{eq:pressurePoisson}  \\
\intertext{with the additional boundary condition,}
  & \grad\cdot\vv =0,  \quad&& \xv\in\partial\OmegaF ,  \label{eq:VPdiv}
\end{alignat}
where the pressure satisfies the Poisson equation~\eqref{eq:pressurePoisson} and the 
boundary condition~\eqref{eq:VPdiv} should be used in addition to the usual conditions
(see~\cite{ICNS} for a discussion of why~\eqref{eq:VPdiv} 
is a suitable {\em pressure} boundary condition).

The evolution of the displacement $\usv=\usv(s,t)$ and velocity $\vsv=\vsv(s,t)$ of the 
structural shell, which depend on arclength $s$ and time $t$, is governed by the equations
\begin{alignat}{3}
  \frac{\partial\usv}{\partial t} &=\vsv, && s\in(0,\len),                    \label{eq:shellEquationU} \\
  \rhos\hs \frac{\partial\vsv}{\partial t} & = \Lsv(\usv) -\sigmav\nv , \qquad&& s\in(0,\len), \label{eq:shellEquation}
\end{alignat}
where the constants $\len$, $\hs$ and $\rhos$ are the length, thickness and mass per unit volume (density) of the shell, respectively, and the
fluid traction, $\sigmav\nv$, is evaluated at $\xv=\xv_0(s)$ on $\GammaI$. 
The vector $\nv$ is the outward unit normal from the fluid domain.
The operator $\Lsv$, in~\eqref{eq:shellEquation} is taken to be 
\begin{align*}
  \Lsv(\usv) = - \Ks\usv+\frac{\partial}{\partial s}\left( \Ts \frac{\partial\usv}{\partial s}\right)
         - \frac{\partial^2}{\partial s^2}\left( \Bs \frac{\partial^2\usv}{\partial s^2}\right),
\end{align*}
where $\Ks$ is a stiffness coefficient, $\Ts$ is a coefficient of tension and $\Bs$ is
a beam coefficient.  More general forms of the shell model are possible~\cite{GuidoboniGlowinskiCavalliniCanic2009}, but the form considered here is sufficient for the purposes of this article.  Initial conditions are required for $\usv(s,0)$ and $\vsv(s,0)$, in general, and suitable boundary conditions at $s=0$ and $s=\len$ are needed as well.

The coupling between the fluid and structure is defined
by the kinematic condition,
\begin{align}
  \vv(\xv_0(s),t)=\vsv(s,t), \quad \xv_0(s)\in \GammaI,\quad s\in(0,\len),   \label{eq:interfaceV}
\end{align}
and the dynamic condition (balance of forces)
which has already been incorporated into the equations through the appearance of the fluid
traction $\sigmav\nv$ in the shell equation~\eqref{eq:shellEquation}.

\section{AMP algorithm}   \label{sec:AMPalgorithm}

In this section we describe our added-mass partitioned (AMP) algorithm to solve
an FSI problem involving an incompressible fluid coupled to a structural shell (or beam).
The essential element of the algorithm is the use of certain AMP interface
conditions, which may be derived at either the continuous or fully discrete
levels.  We start with a discussion of the AMP interface conditions at the
continuous level.  This description leads to different formulations of the
conditions depending on whether the fluid equations are solved numerically
as a fully coupled system or solved using
a fractional-step approach.  In the continuous description, the velocity at the
interface is assigned using a velocity projection.  It is also possible to derive
the AMP interface conditions at a fully discrete level.  An advantage of this
derivation is that the discrete interface conditions maintain the velocity
condition in~\eqref{eq:interfaceV} exactly at a discrete level (without the need
of a projection).  A disadvantage is that the derivation relies on 
a greater knowledge of the discretizations used for
the fluid and shell equations, and so the
approach is less general.

Once the AMP interface conditions are obtained, we then proceed to a description of the AMP algorithm, which is a predictor-corrector method.  The algorithm includes a fractional-step approach in which the solution of the fluid momentum equation in~\eqref{eq:NS3dv} and the pressure equation in~\eqref{eq:pressurePoisson} are updated in separate stages.  For this choice, the AMP interface conditions used in the two separate stages are  identified as the {\em velocity} and {\em pressure} boundary conditions.

The AMP interface conditions are based on the following simple observation. 
\begin{observation}
A Robin boundary condition for the fluid, depending on $\vv$ and $\sigmav$ on the interface $\GammaI$, can be obtained
by substituting the fluid velocity into the shell equation~\eqref{eq:shellEquation} to give 
\begin{align}
&  \rhos\hs \vv_t  = \Lsv(\usv) - \sigmav\nv  , \qquad \xv\in\GammaI.  \label{eq:AMPcondition}
\end{align}
\end{observation}
\noindent 
Equation~\eqref{eq:AMPcondition} can be viewed as generalized Robin (mixed) boundary condition for the fluid of the form $\Bc(\vv,\sigma)=0$,
since the shell displacement $\usv$ is just the time integral of $\vv$ on the interface.

We now consider a basic form of the AMP algorithm that incorporates this boundary condition.  Suppose that the solutions in the fluid and shell are known at a time $t$ and that a partitioned algorithm is desired to advance the solutions to a time $t+\dt$.  In the first step, the equations governing the displacement of the shell are integrated over a time step to obtain $\uvs^{(p)}$, a predicted solution for the shell displacement at $t+\dt$.  Using this predicted displacement in~\eqref{eq:AMPcondition}, we obtain
\begin{equation}
 \sigmav\nv + \frac{\rhos\hs}{\rho}\grad\cdot\sigmav = \Lsv\bigl(\uvs^{(p)}\bigr), \qquad \xv\in\GammaI, 
\label{eq:continuousversion}
\end{equation}
where we have eliminated the fluid acceleration in favor of the divergence of the fluid stress using the momentum equation~\eqref{eq:NS3dv}.  Using the definition of the fluid stress in terms of the pressure and the viscous stress leads to the basic AMP interface condition defined at the continuous level: 
\begin{condition}
The AMP interface condition that can be used when integrating the fluid equations as a fully coupled
system (either the velocity-divergence form~\eqref{eq:NS3dv}--\eqref{eq:NS3dDiv} or 
the velocity-pressure form~\eqref{eq:VPv}--\eqref{eq:pressurePoisson}) is given by
\begin{align}
&   -p\nv  -\frac{\rhos\hs}{\rho}\grad p  + \tauv\nv + \frac{\mu\rhos\hs}{\rho}\Delta \vv  = \Lsv\bigl(\uvs^{(p)}\bigr),  \qquad \xv\in\GammaI.  \label{eq:AMPpressureCondition}
\end{align}
\end{condition}

The interface condition in~\eqref{eq:AMPpressureCondition} is a mixed (Robin) boundary condition, and is the essential component of the AMP algorithm.  The fluid traction term, 
$\sigmav\nv=-p\nv +\tauv\nv$
from the right-hand-side of the structural equations~\eqref{eq:shellEquation} has been explicitly exposed in the fluid 
boundary condition~\eqref{eq:AMPpressureCondition}. The pressure component of this traction term appears as the 
first term on the left of~\eqref{eq:AMPpressureCondition} and is thus effectively being treated in an implicit manner. 

If a fractional-step algorithm is used to integrate the fluid equations, then suitable boundary conditions can be derived by decomposing~\eqref{eq:AMPpressureCondition} into normal and tangential components.  The normal component gives a boundary condition for the pressure, while the tangential component (along with the continuity equation~\eqref{eq:VPdiv}) provide boundary conditions for the velocity:
\begin{condition}
The AMP interface condition that can be used when solving the pressure equation~\eqref{eq:pressurePoisson} is given by
\begin{align}
&   p  + \frac{\rhos\hs}{\rho}\frac{\partial p}{\partial n}  = \nv^T\tauv\nv + \frac{\mu\rhos\hs}{\rho}\nv^T\Delta \vv 
               -\nv^T\Lsv\bigl(\uvs^{(p)}\bigr), \qquad \xv\in\GammaI . \label{eq:AMPpressureBC}
\end{align}
\end{condition}
\begin{condition}
The AMP interface conditions for the velocity that can be used when integrating the momentum equation~\eqref{eq:VPv} are given by
\begin{equation}
\left.
\begin{array}{c}
\displaystyle{
  \tanv^T\tauv\nv + \frac{\mu\rhos\hs}{\rho}\tanv^T\Delta \vv= -\frac{\rhos\hs}{\rho}\tanv^T\grad p  
               + \tanv^T\Lsv\bigl(\uvs^{(p)}\bigr)
}\medskip\\
\displaystyle{
\nabla\cdot\vv=0
}
\end{array}
\right\},  \qquad \xv\in\GammaI,
\label{eq:AMPtangentialVelocityBC}
\end{equation}
where $\tanv$ denotes the unit tangent vector on $\GammaI$ orthogonal to $\nv$. 
\end{condition}

The first condition in~\eqref{eq:AMPtangentialVelocityBC} should be thought of as a boundary 
condition for the tangential component of the velocity, while the continuity equation can be viewed as a boundary condition for the normal component of the velocity. 
The boundary conditions in~\eqref{eq:AMPpressureBC} and~\eqref{eq:AMPtangentialVelocityBC} can be applied
as part of a predictor-corrector algorithm that uses approximate 
values of $\vv$ and $p$ in the right-hand sides of these equations. 
That is, when using~\eqref{eq:AMPpressureBC} to solve for $p$, a guess for $\vv$ is needed, and when using~\eqref{eq:AMPtangentialVelocityBC}
to solve for $\vv$, a guess for $p$ is required.   
It will be shown, through analysis and computations, that the AMP interface conditions can be used to develop
partitioned algorithms that remain stable even when added mass effects are large. 

Having derived the AMP interface conditions at the continuous level, we now consider a derivation at the discrete level for a particular choice of discrete approximations.  Consider the discrete variables,
$\vv^{n}_\iv \approx \vv(\xv_\iv,t^n)$, $p_\iv^n\approx p(\xv_\iv,t^n)$, $\usv^{n}_\iv \approx \usv(s_\iv,t^n)$,
and $\vsv^{n}_\iv \approx \vsv(s_\iv,t^n)$, 
where $\xv_\iv$ denotes the two-dimensional grid points in the fluid domain, and $s_\iv$ denotes the one-dimensional grid points for the shell along $\GammaI$, and $t^n=n\Delta t$.  Here, $\iv=(i_1,i_2)$ is a multi-index for the fluid domain, and where $i_2$ taken to be constant along the interface $\GammaI$.
Consider a particular discrete approximations of the momentum equations for the fluid and shell that uses the trapezodial rule in time, 
\begin{alignat}{3}
& \rhos\hs \frac{\vsv^{n+1}_\iv - \vsv_\iv^n}{\dt} = \frac{1}{2}\Big( \Lsv_{h}(\usv^{n+1}) +  \Lsv_{h}(\usv^{n})\Big)
                  - \frac{1}{2}\Big( \sigmav_\iv^{n+1}\nv_\iv + \sigmav_\iv^{n}\nv_\iv\Big), \qquad&&\iv\in\GammaIh ,\label{eq:solidDiscrete} \\
&  \rho \frac{\vv^{n+1}_\iv - \vv_\iv^n}{\dt} = \frac{1}{2}\Big(\grad_h\cdot\sigmav_\iv^{n+1} + \grad_h\cdot\sigmav_\iv^{n}\Big),
               \quad&& \iv\in\OmegaFh. \label{eq:fluidDiscrete}
\end{alignat}
Here, $\grad_h\cdot$  and $\Lsv_{h}$ denote discrete approximations of the divergence operator and the spatial derivatives in $\Lsv$, 
while $\GammaIh$ and $\OmegaFh$ denote discrete index spaces. 
The approximations given in~\eqref{eq:solidDiscrete} and~\eqref{eq:fluidDiscrete}, which are fully coupled,
can be solved numerically in a partitioned manner by first advancing the solid, using a provisional choice for $\sigmav_\iv^{n+1}\nv_\iv$, and then updating the fluid velocity and pressure,
using a discrete version of the AMP interface condition, which is derived below. 
The solid and
fluid states can be corrected, as needed for accuracy reasons, as part of a predictor-corrector algorithm.

A discrete version of the AMP interface condition may now be derived following
the approach used to obtain the continuous version
in~\eqref{eq:continuousversion}.  Assuming that the discrete fluid and solid
velocities match at $t^n$,
i.e.~$\vv_\iv^n=\vsv_\iv^n$, $\iv\in\GammaIh$, we enforce a similar condition at
$t^{n+1}$ by using it to eliminate the acceleration,
$(\vv_\iv^{n+1}-\vv_\iv^{n})/\dt$, from~\eqref{eq:solidDiscrete} and~\eqref
{eq:fluidDiscrete}, which gives the discrete AMP interface condition
\begin{align}
&  \sigmav_\iv^{n+1}\nv_\iv + \frac{\rhos\hs}{\rho} \grad_h\cdot\sigmav_\iv^{n+1} = 
                \Lsv_{h}\bigl(\usv_\iv^{n+1}\bigr)+  \Lsv_{h}\bigl(\usv_\iv^{n}\bigr) 
      - \sigmav_\iv^{n}\nv_\iv - \frac{\rhos\hs}{\rho}\grad_h\cdot\sigmav_\iv^{n}, \qquad \iv\in\GammaIh.
          \label{eq:pressureBCDiscrete}
\end{align}
As mentioned previously, an advantage of deriving the interface conditions at the discrete level is that 
the discrete fluid and solid velocities match exactly on the interface at each time step assuming they match in the initial conditions.

In general, a more modular approach may be preferred, where the discrete form of the AMP interface conditions do not depend on the choice of the discrete approximations for the equations governing the fluid and shell.  In this case the governing equations can be discretized independently, with appropriate approximations, 
and the AMP interface condition~\eqref{eq:AMPpressureCondition}, or~\eqref{eq:AMPpressureBC} and~\eqref{eq:AMPtangentialVelocityBC} for a fractional-step method,
can be applied when advancing the fluid variables.
In this approach, however, after solving for the new values of $\vv_\iv^{n+1}$ and
$\vvs_\iv^{n+1}$ from the approximations of the fluid and shell equations, it may not be true
that the discrete velocities exactly match on $\GammaIh$.
Therefore, to ensure that the discrete fluid and solid velocities match precisely on $\GammaIh$,
we have found it useful to project the velocities to a common value, based on a
density weighted average, 
\begin{equation}
   \vv^I_\iv = \gamma \vv_\iv^{n+1} + (1-\gamma) \vsv_\iv^{n+1},    \qquad \iv\in\GammaIh, \label{eq:velocityProjection}
\end{equation}
where
\begin{equation}
   \gamma = \frac{1}{1 + (\rhos\hs)/(\rho\hf)} ,
\label{eq:densityWeighting}
\end{equation}
and $\hf$ is a dimensional parameter with units of length that is a measure of the size of the fluid domain.  
In principle,
$\vv_\iv^{n+1}$ and $\vsv_\iv^{n+1}$ are nearly equal and it should not matter
how $\hf$ is chosen.
 In practice, however,
$\vsv_\iv^{n+1}$ may be less reliable when $(\rhos\hs)/(\rho\hf)$ is {\em very
small} due to poor conditioning, while $\vv_\iv^{n+1}$ may be less reliable when
$(\rho\hf)/(\rhos\hs)$ is {\em very small}. 
Therefore
preference is given to the solid velocity when the shell is {\em heavy}, while
giving preference to the fluid velocity when the shell is {\em light}.  In
practice, the results are found to be quite insensitive to the choice of $\hf$.

A second-order accurate AMP algorithm to solve the coupled FSI problem is now described. 
Any number of discrete approximations could be used with the AMP interface conditions. 
The scheme given here is based on the approach
we use for incompressible flow and moving rigid bodies~\cite{splitStep2003,ICNS,mog2006}. 
The fluid equations are solved in velocity-pressure form using a fractional-step algorithm.
Given the discrete solution at the current time $t^n$, and one previous time level, 
the goal is to determine the solution at time $t^{n+1}$. 
A predictor-corrector scheme is used for this purpose as is given in the following algorithm:

\vskip\baselineskip
\noindent{\bf Begin predictor.}

\newcommand{\indentI}{\hangindent\parindent\hangafter=0\noindent}

\medskip
\indentI {\bf Stage I - structure}: Predicted values for the displacement and velocity of the shell are obtained using a second-order accurate {\em leap-frog scheme}
given by
\begin{alignat*}{3}
  \frac{\usv_\iv^{(p)} - \usv_\iv^{n-1}}{2\dt} &= \vsv_\iv^{n}, \qquad&& \iv\in\GammaIh,\\
   \rhos\hs\, \frac{\vsv_\iv^{(p)} - \vsv_\iv^{n-1}}{2\dt} &= \Lsv_{h}\bigl(\usv_\iv^n\bigr) -\sigmav^n_\iv\nv_\iv ,
               \qquad&& \iv\in\GammaIh.
\end{alignat*}

\indentI {\bf Stage II - fluid velocity}: Predicted values for the fluid velocity are obtained using a second-order accurate
 Adams-Bashforth scheme
\begin{alignat*}{3}
  & \rho \frac{\vv_\iv^{(p)} - \vv_\iv^{n}}{\dt} = \frac{3}{2} \Fv_\iv^n - \frac{1}{2} \Fv_\iv^{n-1} , \qquad&& \iv\in\OmegaFh,
\end{alignat*}
where
\begin{alignat*}{3}
  &  \Fv_\iv^n \equiv -\grad_h p_\iv^n + \mu\Delta_h\vv_\iv^n . 
\end{alignat*}
The boundary conditions on $\GammaIh$ are 
\begin{alignat}{3}
 &    \mu\Big( D_{y} v_{1,\iv}^{(p)} + D_{x} v_{2,\iv}^{(p)} \Big) 
   + \frac{\mu\rhos\hs}{\rho} \Delta_{h} v_{1,\iv}^{(p)} = H_\iv,
                 \qquad&& \iv\in\GammaIh, \\
 &  \grad\cdot\vv_\iv^{(p)} =0 ,\quad&& \iv\in\GammaIh,
\end{alignat}
where
\begin{alignat}{3}
 & H_\iv=    -  \frac{\rhos\hs}{\rho}\,\tanv^T\grad p_\iv^{(p)} + \tanv^T\Lsv_{h}\bigl(\usv_\iv^{(p)}\bigr). \label{eq:defH}
\end{alignat}
Here, $D_x$ and $D_y$ are 
some discrete approximations for the first derivatives in the $x$ and $y$ directions, respectively.  A value for
predicted pressure in~\eqref{eq:defH} is not yet available, and so we define it
here using the third-order extrapolation $p_\iv^{(p)} = 3 p_\iv^n - 3
p_\iv^{n-1} + p_\iv^{n-2}$, and then compute a new value in the next stage of
the algorithm.  Appropriate conditions should also be applied on other
boundaries.

\medskip
\indentI{\bf Stage III - fluid pressure}: Predicted values for the pressure are determined by solving 
\begin{align}
   \Delta_h p_\iv^{(p)} & = 0 , \qquad \iv\in\OmegaFh,  \label{eq:pDiscrete}
\end{align}
subject to the boundary conditions, 
\begin{align*}
  p_\iv^{(p)} + \frac{\rhos\hs}{\rho}\, \nv_\iv\cdot\grad_h p_\iv^{(p)} &= 
          \nv_\iv^T\tauv_\iv^{(p)}\nv_\iv  + \frac{\mu\rhos\hs}{\rho}\, \nv_\iv^T\Delta_{h} \vv_\iv^{(p)} - \nv_\iv^T\Lsv_{h}\bigl(\usv_\iv^{(p)}\bigr), 
     \qquad\iv\in\GammaIh ,
\end{align*}
and appropriate conditions applied on the other boundaries.  
Note that in practice it is useful to add a divergence damping term to the right-hand-side of the the pressure 
equation~\eqref{eq:pDiscrete} following~\cite{ICNS,splitStep2003}.

\medskip
\noindent{\bf End predictor.}

\vskip\baselineskip
\noindent{\bf Begin corrector.}

\medskip
\indentI {\bf Stage IV - structure}:
Corrected values for the displacement and velocity of the shell are obtained using a second-order accurate Adams-Moulton scheme, 
\begin{alignat*}{3}
&  \frac{\usv_\iv^{n+1} - \usv_\iv^{n}}{\dt} = \vsv_\iv^{(n+\half)}, \qquad&& \iv\in\GammaIh ,\\
&   \rhos\hs \frac{\vsv_\iv^{n+1} - \vsv_\iv^{n}}{\dt} = \Lsv_{h}\bigl(\usv_\iv^{(n+\half)}\bigr) + \sigmav_\iv^{(n+\half)}\nv_\iv , \qquad&&\iv\in\GammaIh , 
\end{alignat*}
where
\begin{alignat*}{3}
&    \usv_\iv^{(n+\half)} \equiv \frac{ \usv_\iv^{(p)} + \usv_\iv^{n}}{2} ,\qquad \vsv_\iv^{(n+\half)} \equiv \frac{ \vsv_\iv^{(p)} + \vsv_\iv^{n}}{2} ,\qquad 
    ~~ \sigmav_\iv^{(n+\half)} \equiv \frac{ \sigmav_\iv^{(p)} + \sigmav_\iv^{n}}{2}.  
\end{alignat*}

\indentI {\bf Stage V - fluid velocity}: Corrected values for the fluid velocity are obtained using a second-order accurate Adams-Moulton scheme,
\begin{align*}
   \rho \frac{\vv_\iv^{n+1} - \vv_\iv^{n}}{\dt} &= \frac{1}{2} ( \Fv_\iv^{(p)} + \Fv_\iv^{n} ), \qquad \iv\in\OmegaFh.
\end{align*}
The boundary conditions have the same form as those used in Stage II with $\vv_\iv^{(p)}$ and $\usv_\iv^{(p)}$ replaced by $\vv_\iv^{n+1}$ and $\usv_\iv^{n+1}$.

\medskip
\indentI {\bf Stage VI - fluid pressure}: Corrected values for the pressure are determined by solving the discrete equations in Stage III with the predicted values replaced by values at $t^{n+1}$.

\medskip
\indentI {\bf Stage VII - project interface velocity}: Perform a projection so that the fluid and solid velocities match on the interface, using
\begin{align*}
   \vv_\iv^I = \gamma \vv_\iv^{n+1} + (1-\gamma) \vsv_\iv^{n+1}, \qquad \iv\in\GammaIh,
\end{align*}
where $\gamma$ is defined in~\eqref{eq:densityWeighting}, and then set
\begin{align*}
   \vsv_\iv^{n+1} = \vv_\iv^{n+1} = \vv_\iv^I,\qquad \iv\in\GammaIh.
\end{align*}
\medskip
\noindent{\bf End corrector.}

\medskip
We emphasize that the AMP algorithm is stable with no corrector step, although if the predictor step is used alone, then Stage VII should be performed
following the predictor to project the interface velocity. 
We typically use the corrector step, since for the fluid in isolation the 
scheme has a larger stability region than the predictor alone and
the stability region includes the imaginary axis so that the scheme can be used for inviscid problems ($\mu=0$).

\section{FSI Model problems} \label{sec:modelProblems}

In this section, three FSI model problems of increasing complexity are defined. 
The simplest model problem, MP-I1, involves an inviscid incompressible fluid and a shell model that only supports vertical motion.  The second and third models involve a viscous compressible fluid coupled to a shell that supports only vertical motion, MP-V1, and motion in both the horizontal and vertical directions, MP-V2.  The first model, MP-I1, is used to study the stability of the
AMP scheme in Section~\ref{eq:analysis}.
Exact traveling wave solutions to all three model problems are
given in~\ref{sec:travelingWave}, and numerical calculations of the model problems are
given in Section~\ref{sec:results}. 
In all cases the fluid domain is the rectangle $\OmegaF=(0,L)\times(-H,0)$. The shell is defined on the interval $\OmegaS=\{ x \,\vert\, x\in(0,L)\}$, and it is parameterized by $x$ (which is equal to the arclength $s$). The solution is assumed to be periodic in the $x$-direction
with period $L$. The domain geometry is shown in Figure~\ref{fig:elasticShellCartoon}.

%
%
{
\newcommand{\yH}{3}
\newcommand{\xL}{8}
\begin{figure}[hbt]
\newcommand{\textFont}{\normalss}
\begin{center}
\begin{tikzpicture}[scale=1]
\useasboundingbox (0,0) rectangle (\xL,3.5);  
\draw[-,thick] (0,\yH) -- (\xL,\yH);
\draw[-,thick] (\xL,0) -- (\xL,\yH);
\draw[-,thick] (0,0) node[anchor=north] {$x=0$} -- (\xL,0) node[anchor=north] {$x=L$};
\draw[-,thick] (0,0) node[anchor=east] {$y=-H$} -- (0,\yH) node[anchor=east] {$y=0$};
\draw[thick,blue,yshift=\yH cm] (0,0) sin (2,0.25) cos (4,0) sin (6,-.25) cos (8,0); 
\draw[blue] (4,\yH) node[anchor=south,yshift=5pt] {\textFont $y=\eta(x,t)$};
%
\end{tikzpicture}
\end{center}
\caption{The geometry for the 2D fluid domain and structural shell or beam.}
\label{fig:elasticShellCartoon}
\end{figure}
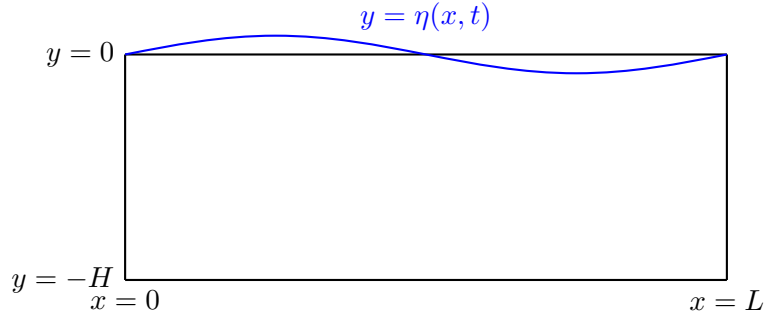
}

\begin{ModelProblemI}
The first model problem defines an inviscid incompressible fluid
and a shell that only supports vertical motion (with $\eta=\us_2$), 
\begin{equation*}
\begin{aligned}
& \text{Fluid:} \quad \left\{ 
   \begin{alignedat}{3}
  &  \rho\frac{\partial\vv}{\partial t} 
                 + \grad p =0, \qquad&& \xv\in\OmegaF ,  \\
  & \grad\cdot\vv =0,  \qquad&& \xv\in\OmegaF ,   \\
  & v_2(x,-H,t)=0 , ~~&&\text{or}~~ p(x,-H,t)=0  , \qquad x\in(0,L), 
   \end{alignedat}  \right. 
\\ 
& \text{Shell:}\quad  \left\{ 
   \begin{alignedat}{3}
  &  \rhos\hs\frac{\partial^2\eta}{\partial t^2}  = -\Ks \eta + \Ts\frac{\partial^2\eta}{\partial x^2} + p(x,0,t), \qquad&& x\in\OmegaS,
   \end{alignedat}  \right. \\
& \text{Interface:}\quad   v_2(x,0,t)=\frac{\partial\eta}{\partial t}(x,t), \qquad x\in\OmegaS .
\end{aligned}
  \label{eq:MP1}
\end{equation*}
\end{ModelProblemI}

\begin{ModelProblemVa}
The second model problem defines a viscous incompressible fluid 
and a shell that only supports vertical motion (with $\eta=\us_2$),
\begin{equation*}
\begin{aligned}
& \text{Fluid:}\quad  \left\{ 
   \begin{alignedat}{3}
  &  \rho \frac{\partial\vv}{\partial t} 
                  +  \grad p = \mu\Delta\vv, \qquad&& \xv\in\OmegaF ,  \\
  & \grad\cdot\vv =0,  \qquad&& \xv\in\OmegaF , \\
  & \vv(x,-H,t)=0 ,  \qquad&& x\in(0,L), 
   \end{alignedat}  \right. 
\\
& \text{Shell:}\quad \left\{ 
   \begin{alignedat}{3}
  &  \rhos\hs\frac{\partial^2\eta}{\partial t^2}  = -\Ks \eta + \Ts\frac{\partial^2\eta}{\partial x^2}
                 + p(x,0,t) -2\mu\frac{\partial v_2}{\partial y}(x,0,t), \qquad&& x\in\OmegaS, 
   \end{alignedat}  \right.  
 \\
& \text{Interface:}\quad v_1(x,0,t)=0, \qquad v_2(x,0,t)=\frac{\partial\eta}{\partial t}(x,t), \qquad x\in\OmegaS .
\end{aligned}
  \label{eq:MP2}
\end{equation*}
\end{ModelProblemVa}

\begin{ModelProblemVb}
The third model problem defines a viscous incompressible fluid
and a shell that supports horizontal and vertical motion,
\begin{equation*}
\begin{aligned}
& \text{Fluid:}  \quad
  \left\{ 
   \begin{alignedat}{3}
  &  \rho \frac{\partial\vv}{\partial t} 
                 + \grad p = \mu\Delta\vv, \qquad&& \xv\in\OmegaF ,  \\
  & \grad\cdot\vv =0,  \qquad&& \xv\in\OmegaF  , \\
  & \vv(x,-H,t)=0 ,  \qquad&& x\in(0,L),   
   \end{alignedat}  \right. 
\\ 
& \text{Shell:} \quad    \left\{ \displaystyle
   \begin{alignedat}{3}
  &  \rhos\hs\frac{\partial^2\usv}{\partial t^2}  = -\Ks \usv + \Ts\frac{\partial^2\usv}{\partial x^2}
                 + \begin{bmatrix} -\mu( \partial_y v_1 + \partial_x v_2 ) \\
                                   p -2\mu \partial_y v_2 \end{bmatrix} (x,0,t)
                 , \qquad&& x\in\OmegaS, 
   \end{alignedat}  \right.  
 \\
& \text{Interface:}  \quad \vv(x,0,t)=\frac{\partial\usv}{\partial t}(x,t) , \qquad x\in\OmegaS.
\end{aligned}
  \label{eq:MP3}
\end{equation*}
\end{ModelProblemVb}

\section{Analysis of an inviscid fluid and simplified shell} \label{eq:analysis}

In this section, a mode analysis is used to study the stability of a traditional
partitioned scheme and the new added-mass partitioned scheme for the model
problem MP-I1.  For both schemes, we advance the fluid equations using a
backward-Euler scheme as opposed to the predictor-corrector time-stepping scheme
used in the AMP algorithm described in Section~\ref{sec:AMPalgorithm}.  While
this simplifies the analysis somewhat, the essential results are unaffected.  As
expected, we find that the traditional scheme becomes unconditionally unstable
for light structures.  The AMP scheme, in contrast, remains stable for
arbitrarily light structures.
The traditional partitioned scheme has been analysed previously by  Causin, Gerbeau and Nobile~\cite{CausinGerbeauNobile2005}, among others.
The results here extend that work. 

To begin the analysis, the solution is expanded in a finite Fourier series in $x$ (i.e.~a pseudo-spectral approximation), which gives 
\[
 \vv(x,y,t) \approx \sum_{k=-M}^{M} \vvh(k,y,t) e^{2\pi i k x/L} , \qquad
    p(x,y,t) \approx \sum_{k=-M}^{M} \ph(k,y,t) e^{2\pi i k x/L} ,
\]
and
\[
  \eta(x,t) \approx \sum_{k=-M}^{M} \etah(k,t) e^{2\pi i k x/L} . 
\]
The governing equations for the Fourier components corresponding to a given wavenumber $k$ are
\[
\left. \begin{array}{r}
\displaystyle{
  \rho \partial_t \vh_1 + i\kx \ph =0
} \smallskip\\
\displaystyle{
  \rho \partial_t \vh_2 + \ph_y =0
} \smallskip\\
\displaystyle{
  \ph_{yy} - \kx^2 \ph =0
}
\end{array}
\right\},\qquad y \in (-\Hf,0),
\]
for the fluid, and
\begin{equation}
  \rhos\hs \etah_{tt} = - \Lt \etah + \ph(k,0,t) ,
\label{eq:modelShellEquationFT}
\end{equation}
for the shell, where $\vvh=(\vh_1,\vh_2)$, $\kx = 2\pi k/L$ and $\Lt=K_s+T_s\kx^2+B_s\kx^4$.
The boundary condition on the fluid at $y=-\Hf$ is taken as a solid {\em slip} wall (or symmetry condition), which implies $\vh_2(k,-\Hf,t)=0$.  It then follows from the fluid momentum equation that $\ph_y(k,-\Hf,t)=0$.
The conditions used at the interface depend on the choice of algorithm, as described below, and in principle a full solution would require a choice for the initial conditions.


We discretize in time but for clarity leave the spatial coordinate $y$ continuous 
(the analysis can be easily extended to fully discrete in $y$). 
Let
\[
\vvh^n(k,y)\approx\vvh(k,y,t^n), \qquad \ph^n(k,y)\approx\ph(k,y,t^n), \qquad \etah^n(k)\approx\etah(k,t^n),
\]
 where $t^n=n\dt$. 
Let $D_{+t}$ and $D_{-t}$ denote the forward and backward divided different operators in time,
e.g.~$D_{+t}\etah^n=(\etah^{n+1}-\etah^n)/\dt$ and $D_{-t}\etah^n=(\etah^{n}-\etah^{n-1})/\dt$.

We begin with an analysis of the traditional coupling scheme, which is described in the following algorithm:

\begin{algorithm} {\bf Traditional partitioned scheme:} 

\noindent Stage I: Advance the shell displacement using a leap-frog scheme:
\begin{align}
  \rhos\hs D_{+t}D_{-t}\etah^n & = -\Lt \etah^n + \ph^n(k,0) . \label{eq:etaLeapFrog}
\end{align}
Stage II: Advance the fluid velocity and pressure using a backward-Euler scheme,
\begin{equation}
\left. \begin{array}{r}
\displaystyle{
\rho \frac{\vh_1^{n+1}-\vh_1^{n}}{\dt}  + i\kx \ph^{n+1} =0
} \smallskip\\
\displaystyle{
\rho \frac{\vh_2^{n+1}-\vh_2^{n}}{\dt} + \ph_y^{n+1} =0
} \smallskip\\
\displaystyle{
\ph_{yy}^{n+1} - \kx^2 \ph^{n+1} =0
}
\end{array} \right\}, \quad y \in (-\Hf,0),
\label{eq:fhat}
\end{equation}
with interface and boundary conditions,
\begin{equation}
\begin{array}{ll}
\ph_y^{n+1} = -\rho  D_{+t}D_{-t}\etah^n , \qquad & y=0, \smallskip\\
\ph_y^{n+1}=0, \qquad &  y=-\Hf.
\end{array}
\label{eq:pressureBCTraditional}
\end{equation}
\label{alg:TP}
\end{algorithm}

For $\kx\ne 0$, the solution of the equation in~\eqref{eq:fhat} for $\ph^{n+1}$, with interface and boundary conditions in~\eqref{eq:pressureBCTraditional}, is given by
\begin{align}
  \ph^{n+1} &= -\rho D_{+t}D_{-t}\etah^n \frac{\cosh(\kx (y+\Hf))}{\kx \sinh(\kx \Hf)}. \label{eq:IESpressure}
\end{align}
For $\kx=0$, assume $\etah^n(0)=0$ and
 $\ph^n(0,y)=0$ (so that a flat interface at $y=0$ is in equilibrium when the pressure is zero). 
From~\eqref{eq:IESpressure}, the pressure on the interface, $y=0$, is 
\begin{align}
  \ph^{n+1}(k,0) &= - M_a D_{+t}D_{-t}\etah^n, \label{eq:pTraditional}
\end{align}
where $M_a$ is an added-mass coefficient given by, 
\begin{align}
   M_a=\frac{\rho\cosh(\kx \Hf)}{\kx \sinh(\kx \Hf)} 
      =\frac{\rho L}{2\pi k} \, \coth(2\pi k \Hf/L). \label{eq:Ma}
\end{align}
Substituting~\eqref{eq:pTraditional} into~\eqref{eq:etaLeapFrog} 
gives a difference equation for $\etah^{n}$, 
\begin{align}
  \rhos\hs D_{+t}D_{-t}\etah^n & = -\Lt \etah^n  
       - M_a ~D_{+t}D_{-t}\etah^{n-1}.  \label{eq:etaTraditional} 
\end{align}
The stability of the overall scheme is determined by the stability of~\eqref{eq:etaTraditional}. 
The appearance of the factor involving $M_a$ on the right-hand-side of~\eqref{eq:etaTraditional} (i.e. lagged in time) is the primary
source of the added-mass instability.

Assuming solutions of the form $\etah^{n}=\amp^n \etah^{0}$, it follows that the amplification
factor $\amp$ satisfies the cubic equation
\begin{align}
  \frac{\rhos\hs}{\dt^2} (\amp-1)^2 \amp + \Lt \amp^2 + \frac{M_a}{\dt^2} (\amp-1)^2 &=0. \label{eq:FSIpolytraditional}
\end{align}
The scheme is called weakly stable if all roots to~\eqref{eq:FSIpolytraditional} satisfy $\vert\amp\vert\le 1$. 
Using the theory of von Neumann polynomials~\cite{Miller1971,Strikwerda89} one can show the following necessary and
sufficient condition.
\begin{theorem} \label{th:traditional}
  Assuming $\rhos\hs>0$, $\Lt>0$ and $M_a>0$, 
the traditional partitioned scheme given in Algorithm~\ref{alg:TP} is weakly stable if and only if 
the following conditions are satisfied:
\begin{align}
  &  M_a < \rhos \hs, \label{eq:traditionalStabLimitI} \\
  &  \dt^2 <  4 \, \frac{\rhos\hs}{\Lt} \Big[ 1- \frac{M_a}{\rhos\hs} \Big]. \label{eq:traditionalStabLimitII}
\end{align}
\end{theorem}
The condition~\eqref{eq:traditionalStabLimitI} agrees with the result of Causin, Gerbeau and Nobile~\cite{CausinGerbeauNobile2005}.
The additional condition~\eqref{eq:traditionalStabLimitII} indicates how the time-step must be reduced as the added-mass effects become larger.
It is noted that the traditional partitioned scheme is unconditionally unstable if the contribution from the added mass is too large, i.e.~if $M_a >\rhos \hs$, regardless of the choice for the time-step, $\dt$.  From~\eqref{eq:Ma}, we observe that the added-mass coefficient is large for large $\rho L$, small $k$ and/or small $\Hf/L$.  Thus, the traditional partitioned scheme has difficulty when the mass ratio, $(\rho L)/(\rhos\hs)$, is large, the wave number of surface variations is low, and/or the fluid domain is thin.

We now contrast the stability results found for the traditional partitioned scheme with the corresponding ones for the new {\em added-mass} partitioned (AMP) algorithm.  This is done by considering the following  version of the AMP algorithm:

\begin{algorithm} {\bf AMP scheme:} 

\noindent Stage I: Advance the shell displacement using a leap-frog scheme (as in the traditional scheme),
\begin{equation}
  \rhos\hs D_{+t}D_{-t}\etah^n  = -\Lt \etah^n + \ph^n(k,0) .
\label{eq:etaLeapFrogAM}
\end{equation}
Stage II: Advance the fluid velocity and pressure using a backward-Euler scheme (as in the traditional scheme),
\[
\left. \begin{array}{r}
\displaystyle{
\rho \frac{\vh_1^{n+1}-\vh_1^{n}}{\dt}  + i\kx \ph^{n+1} =0
} \smallskip\\
\displaystyle{
\rho \frac{\vh_2^{n+1}-\vh_2^{n}}{\dt} + \ph_y^{n+1} =0
} \smallskip\\
\displaystyle{
\ph_{yy}^{n+1} - \kx^2 \ph^{n+1} =0
}
\end{array} \right\}, \quad y \in (-\Hf,0),
\]
with the AMP Robin interface condition, derived from~\eqref{eq:AMPpressureBC},
\begin{equation}
\ph^{n+1} + \frac{\rhos\hs}{\rho} \ph_y^{n+1} = \Lt \etah^{n+1} , \qquad y=0,
\label{eq:AMpressureBC}
\end{equation}
and boundary conditions,
\[
\ph_y^{n+1}=0, \qquad  y=-\Hf.
\]
\label{alg:AMP}
\end{algorithm}

The key ingredient of the AMP scheme is the Robin condition in~\eqref{eq:AMpressureBC} for the fluid pressure.  The solution for the pressure is now given by
\begin{align}
   \ph^{n+1} &= \Lt\etah^{n+1}\frac{\rho\cosh(\kx (y+\Hf))}{\rho\cosh(\kx \Hf) + \rhos\hs \kx \sinh(\kx \Hf)}. \label{eq:pressureAM}
\end{align}
Combining~\eqref{eq:etaLeapFrogAM} and~\eqref{eq:pressureAM} gives the following difference equation for $\etah^n$:
\begin{equation}
  \rhos\hs D_{+t}D_{-t}\etah^n  = -\Lt \etah^n + \Lt\etah^{n}\frac{\rho\cosh(\kx \Hf)}{\rho\cosh(\kx \Hf) + \rhos\hs \kx \sinh(\kx \Hf)} = -\Lt B \etah^n
\label{eq:etaAM}
\end{equation}
where 
\begin{align*}
   B &\equiv \frac{\rhos\hs \kx \sinh(\kx \Hf)}{ \rho \cosh(\kx \Hf) + \rhos\hs \kx \sinh(\kx \Hf) }
     =  \frac{1}{1+M_a}   .
\end{align*}
Note that $0< B< 1$ and thus solutions to~\eqref{eq:etaAM} have a reduced frequency in time, $\sqrt{\Lt B/(\rhos\hs)}$, as
compared to the natural frequency, $\sqrt{\Lt/(\rhos\hs)}$, of the shell with no fluid present. 

\begin{theorem}
   The added-mass partitioned scheme described in Algorithm~\ref{alg:AMP} is stable if and only if 
\begin{align}
   \dt < 2 \sqrt{\frac{\rhos\hs}{\Lt B}} 
           = 2 \sqrt{\frac{\rhos\hs(1+M_a)}{\Lt}} . \label{eq:AMPstability}
\end{align}
\end{theorem}
\begin{proof}
Using $\etah^{n}=\amp^n \etah^{0}$ and substituting into~\eqref{eq:etaAM} gives
\begin{align*}
&   \amp^2 -2(1-b\dt^2) \amp + 1 = 0, 
\end{align*}
where  $b = (\Lt B)/(2\rhos\hs)$.
The requirement that $\vert\amp\vert \le 1$ and that the roots be simple 
implies $\vert 1-b\dt^2\vert < 1$ which leads to~\eqref{eq:AMPstability}.
\end{proof}

It is interesting to note that the AMP algorithm (for the MP-I1 model problem) is stable when using the time-step restriction that arises when advancing the shell with no fluid present (i.e.~taking $M_a=0$ in~\eqref{eq:AMPstability}).  Also, when~\eqref{eq:AMPstability} holds the amplification factor satisfies $\vert\amp\vert=1$ and the scheme is non-dissipative.  Finally, for this simple AMP scheme, the fluid and shell velocities 
at the interface only match to first-order accuracy. As noted previously, one
could include an additional projection step to match the velocities on the interface.


\renewcommand{\tableFontSize}{\scriptsize}
\section{Numerical results} \label{sec:results}

Numerical results are now presented that verify the accuracy and stability of the AMP algorithm.
Solutions to the three model problems defined in Section~\ref{sec:modelProblems}
are computed for exact solutions constructed with the method of analytic solutions, 
and for exact traveling wave solutions provided in~\ref{sec:travelingWave}.
The numerical scheme is a predictor-corrector algorithm that follows the
algorithm given in Section~\ref{sec:AMPalgorithm}.
Note that the AMP algorithm remains stable with no corrector step provided 
the predictor in the fluid is stable in isolation (the Adams-Bashforth predictor requires
$\mu>0$ or sufficient dissipation).

In all cases the parameters are chosen as, 
\begin{align*}
  \rho=1,\quad H=1, \quad \rhos\hs=\Ts = \dsf\,\rho H, \quad \Ks=0, \quad\Bs=0,
\end{align*}
where the {\em density ratio} $\dsf=\rhos\hs/(\rho H)$ is varied to represent both light and heavy structures.
The fluid length scale in the velocity projection~\eqref{eq:densityWeighting} is taken as $\hf=10$, although
the results are very insensitive to this choice. 
The grid for the fluid domain is a Cartesian grid for the rectangular region $\OmegaF=(0,L)\times(-H,0)$
with $N_j+1$ grid points in each direction and grid spacing $h_j=\dx_j=L/N_j =\dy_j=H/N_j$. 
The one-dimensional grid
for the shell has $N_j+1$ grid points. 
Here the 
subscript $j$ will be used to note grids of varying resolution.


%
\subsection{The method of analytic solutions}

\bogus{
\begin{table}[hbt]\tableFont 
\begin{center}
\begin{tabular}{|l|c|c|c|c|c|c|c|c|c|} \hline 
grid  & N & $p$ & r &  $\vert \mathbf{v} \vert$ & r & $\bar{\mathbf{u}}$ & r & $\bar{\mathbf{v}}$ & r \\ \hline 
                sq20 &     1 & \num{6.0}{-2} &      & \num{2.7}{-2} &      & \num{2.0}{-2} &      & \num{2.7}{-2} &       \\ \hline
                sq40 &     2 & \num{1.3}{-2} &  4.5 & \num{6.0}{-3} &  4.4 & \num{5.0}{-3} &  4.0 & \num{6.0}{-3} &  4.4  \\ \hline
                sq80 &     4 & \num{3.0}{-3} &  4.4 & \num{1.4}{-3} &  4.4 & \num{1.2}{-3} &  4.1 & \num{1.4}{-3} &  4.4  \\ \hline
               sq160 &     8 & \num{7.3}{-4} &  4.1 & \num{3.2}{-4} &  4.2 & \num{3.0}{-4} &  4.1 & \num{3.2}{-4} &  4.2  \\ \hline
    rate             &       &  $2.12$       &      &  $2.12$       &      &  $2.02$       &      &  $2.12$       &       \\ \hline
\end{tabular}
\caption{Inses shell: fw.mup05.full.am1.rep2, AM=1, amImp=1, vpf=0.1, ts=pc2, $t=.5$, $\mu=0.05$, $\rho=1$, $\bar{\rho}=1e2$, T=1e2, , cfl=0.75, normalMotion=0, gp=pn, Mon Jul 29  8:14:15 2013}\label{table:shell.fw.mup05.full.am1.rep2.pc2.pc2}
\end{center}
\end{table}
}
\newcommand{\tableTZMPVIIAmpHeavy}{%
\begin{tabular}{|c|c|c|c|c|c|c|c|c|} \hline
\multicolumn{9}{|c|}{MP-V2, trigonometric solution, heavy solid} \\ \hline 
\strutt~~$h_j$~~& $E_j^{(p)}$ & $r$ & $E_j^{(\vv)}$ & $r$ & $E_j^{(\usv)}$ & $r$  & $E_j^{(\vsv)}$  & $r$  \\ \hline 
 1/20      &  \num{6.0}{-2} &      & \num{2.7}{-2} &      & \num{2.0}{-2} &      & \num{2.7}{-2} &       \\ \hline
 1/40      &  \num{1.3}{-2} &  4.5 & \num{6.0}{-3} &  4.4 & \num{5.0}{-3} &  4.0 & \num{6.0}{-3} &  4.4  \\ \hline
 1/80      &  \num{3.0}{-3} &  4.4 & \num{1.4}{-3} &  4.4 & \num{1.2}{-3} &  4.1 & \num{1.4}{-3} &  4.4  \\ \hline
 1/160	   &  \num{7.3}{-4} &  4.1 & \num{3.2}{-4} &  4.2 & \num{3.0}{-4} &  4.1 & \num{3.2}{-4} &  4.2  \\ \hline
\rateLabel &   $2.12$       &      &  $2.12$       &      &  $2.02$       &      &  $2.12$       &       \\ \hline
\end{tabular}
}

\bogus{
\begin{table}[hbt]\tableFont 
\begin{center}
\begin{tabular}{|l|c|c|c|c|c|c|c|c|c|} \hline 
grid  & N & $p$ & r &  $\vert \mathbf{v} \vert$ & r & $\bar{\mathbf{u}}$ & r & $\bar{\mathbf{v}}$ & r \\ \hline 
                sq20 &     1 & \num{1.2}{-2} &      & \num{6.0}{-3} &      & \num{7.1}{-3} &      & \num{7.1}{-3} &       \\ \hline
                sq40 &     2 & \num{3.7}{-3} &  3.4 & \num{1.3}{-3} &  4.5 & \num{1.7}{-3} &  4.2 & \num{1.7}{-3} &  4.2  \\ \hline
                sq80 &     4 & \num{1.0}{-3} &  3.6 & \num{3.1}{-4} &  4.2 & \num{4.2}{-4} &  4.0 & \num{4.2}{-4} &  4.0  \\ \hline
               sq160 &     8 & \num{2.6}{-4} &  3.9 & \num{7.8}{-5} &  4.0 & \num{1.1}{-4} &  4.0 & \num{1.1}{-4} &  4.0  \\ \hline
    rate             &       &  $1.86$       &      &  $2.08$       &      &  $2.03$       &      &  $2.03$       &       \\ \hline
\end{tabular}
\caption{Inses shell: tz.trig.am1.rep0, AM=1, amImp=1, vpf=0.1, ts=pc2, $t=.5$, $\mu=0.05$, $\rho=1$, $\bar{\rho}=1$, T=1, TZ=trig(fx=2.,ft=2), cfl=0.75, normalMotion=0, gp=pn, Mon Jul 29  9:55:30 2013}\label{table:shell.tz.trig.am1.rep0.pc2.pc2}
\end{center}
\end{table}
}
\newcommand{\tableTZMPVIIAmpMedium}{%
\begin{tabular}{|c|c|c|c|c|c|c|c|c|} \hline
\multicolumn{9}{|c|}{MP-V2, trigonometric solution, medium solid} \\ \hline 
\strutt~~$h_j$~~& $E_j^{(p)}$ & $r$ & $E_j^{(\vv)}$ & $r$ & $E_j^{(\usv)}$ & $r$  & $E_j^{(\vsv)}$  & $r$  \\ \hline 
 1/20      &  \num{1.2}{-2} &      & \num{6.0}{-3} &      & \num{7.1}{-3} &      & \num{7.1}{-3} &       \\ \hline
 1/40      &  \num{3.7}{-3} &  3.4 & \num{1.3}{-3} &  4.5 & \num{1.7}{-3} &  4.2 & \num{1.7}{-3} &  4.2  \\ \hline
 1/80      &  \num{1.0}{-3} &  3.6 & \num{3.1}{-4} &  4.2 & \num{4.2}{-4} &  4.0 & \num{4.2}{-4} &  4.0  \\ \hline
 1/160	   &  \num{2.6}{-4} &  3.9 & \num{7.8}{-5} &  4.0 & \num{1.1}{-4} &  4.0 & \num{1.1}{-4} &  4.0  \\ \hline
\rateLabel &   $1.86$       &      &  $2.08$       &      &  $2.03$       &      &  $2.03$       &       \\ \hline
\end{tabular}
}

\bogus{
\begin{table}[hbt]\tableFont 
\begin{center}
\begin{tabular}{|l|c|c|c|c|c|c|c|c|c|} \hline 
grid  & N & $p$ & r &  $\vert \mathbf{v} \vert$ & r & $\bar{\mathbf{u}}$ & r & $\bar{\mathbf{v}}$ & r \\ \hline 
                sq20 &     1 & \num{1.2}{-2} &      & \num{3.9}{-3} &      & \num{5.3}{-3} &      & \num{5.3}{-3} &       \\ \hline
                sq40 &     2 & \num{3.7}{-3} &  3.1 & \num{8.1}{-4} &  4.8 & \num{8.0}{-4} &  6.6 & \num{8.1}{-4} &  6.6  \\ \hline
                sq80 &     4 & \num{1.0}{-3} &  3.7 & \num{1.9}{-4} &  4.3 & \num{1.7}{-4} &  4.7 & \num{1.9}{-4} &  4.3  \\ \hline
               sq160 &     8 & \num{2.6}{-4} &  3.9 & \num{4.8}{-5} &  4.0 & \num{4.2}{-5} &  4.1 & \num{4.8}{-5} &  4.0  \\ \hline
    rate             &       &  $1.84$       &      &  $2.12$       &      &  $2.32$       &      &  $2.25$       &       \\ \hline
\end{tabular}
\caption{Inses shell: tz.trig.am1.rem2, AM=1, amImp=1, vpf=0.1, ts=pc2, $t=.5$, $\mu=0.05$, $\rho=1$, $\bar{\rho}=.01$, T=.01, TZ=trig(fx=2.,ft=2), cfl=0.75, normalMotion=0, gp=pn, Mon Jul 29 10:10:10 2013}\label{table:shell.tz.trig.am1.rem2.pc2.pc2}
\end{center}
\end{table}

}
\newcommand{\tableTZMPVIIAmpLight}{%
\begin{tabular}{|c|c|c|c|c|c|c|c|c|} \hline
\multicolumn{9}{|c|}{MP-V2, trigonometric solution, light solid} \\ \hline 
\strutt~~$h_j$~~& $E_j^{(p)}$ & $r$ & $E_j^{(\vv)}$ & $r$ & $E_j^{(\usv)}$ & $r$  & $E_j^{(\vsv)}$  & $r$  \\ \hline 
 1/20      &  \num{1.2}{-2} &      & \num{3.9}{-3} &      & \num{5.3}{-3} &      & \num{5.3}{-3} &       \\ \hline
 1/40      &  \num{3.7}{-3} &  3.1 & \num{8.1}{-4} &  4.8 & \num{8.0}{-4} &  6.6 & \num{8.1}{-4} &  6.6  \\ \hline
 1/80      &  \num{1.0}{-3} &  3.7 & \num{1.9}{-4} &  4.3 & \num{1.7}{-4} &  4.7 & \num{1.9}{-4} &  4.3  \\ \hline
 1/160	   &  \num{2.6}{-4} &  3.9 & \num{4.8}{-5} &  4.0 & \num{4.2}{-5} &  4.1 & \num{4.8}{-5} &  4.0  \\ \hline
\rateLabel &   $1.84$       &      &  $2.12$       &      &  $2.32$       &      &  $2.25$       &       \\ \hline
\end{tabular}
}

\begin{figure}[hbt]\tableFontSize
\begin{center}
\tableTZMPVIIAmpHeavy
\vskip.5\baselineskip
\tableTZMPVIIAmpMedium
\vskip.5\baselineskip
\tableTZMPVIIAmpLight
\caption{Trigonometric exact solution for a viscous incompressible fluid and structural shell (model problem MP-V2).
Maximum errors and estimated convergence rates at $t=0.5$ computed using the AMP scheme for a heavy solid, $\delta=10^3$,
medium solid, $\delta=1$, and light solid, $\delta=10^{-2}$, where $\delta=(\rhos\hs)/(\rho H)$.
}
\label{tab:TZ_MP-V2_AMP}
\end{center}
\end{figure}

The method of analytic solutions is a useful technique for constructing exact
solutions of initial-boundary-value problems for partial differential equations
for the purpose of checking the behavior and accuracy of the numerical
implementation of a problem.  This method, also known as the {\it method of
manufactured solutions}~\cite{Roache2002} or {\it twilight-zone
forcing}~\cite{CGNS}, adds forcing functions to the governing equations and
boundary conditions. These forcing functions are specified so that a chosen
function, $\tilde{\qv}(\xv,t)$, becomes the exact solution of the forced equations,
and thus the error in the discrete solution can be computed exactly.

The method of analytic solutions is applied to our FSI problem using
trigonometric functions for the components of  fluid velocity and pressure given by
\begin{alignat}{4}
\vtz_1 &= \half\cos(\fx\pi x)\cos(\fx\pi y)\cos(\ft\pi y), \\
\vtz_1 &= \half\sin(\fx\pi x)\sin(\fx\pi y)\cos(\ft\pi y), \\
\ptz   &=      \cos(\fx\pi x)\cos(\fx\pi y)\cos(\ft\pi y),
\end{alignat}
where $\fx$, $\ft$ are frequency parameters.
The velocity solution is chosen to be divergence free.
The analytic solution for the simplified shell is chosen to be a standing wave
\begin{align}
  \ustz_1 &=   \bar{a}\sin( \fx \pi x ) \cos(f_x \pi c t ) , \\
  \ustz_2 &=   \bar{b}\sin( \fx \pi x ) \cos(f_x \pi c t ) , 
\end{align}
with frequency $f_x$, amplitude $\bar{a}=\bar{b}=0.1$ and speed $c=(\Ts/\rhos)^{1/2}$ .
These forms are substituted into the governing equations and boundary conditions to define the
required forcing functions.

Numerical solutions of the (forced) model problems are computed using the AMP
algorithm.  The initial conditions for the numerical calculations are taken from
the exact solutions. 
For each case, maximum-norm errors,
$E_j^{(q)}$, for solution component $q$, are computed on grids of increasing
resolution using grid spacings $\dx_j=\dy_j=h_j=1/(20j)$, $j=1,2,\ldots$.  The
convergence rate, $\zeta_q$, is estimated by a least squares fit to the
logarithm of the error equation, $E_j^{(q)} = C_q h_j^{\zeta_q}$, where $C_q$ is
approximately constant for small grid spacings.  For vector variables, such as
$\vv$ or $\usv$, the error denotes the maximum over all components of the
vector.

The tables in Figure~\ref{tab:TZ_MP-V2_AMP} show the maximum-norm error and
estimated convergence rates at $t=0.5$ for the density ratios $\dsf=10^{-2}$
(light shell), $\dsf=1$ (medium shell) and $\dsf=10^{3}$ (heavy shell), using
the parameters $\fx=\ft=2$ and $\mu=0.05$.  The values in the columns labeled
``r'' give the ratio of the error on the current grid to that on the previous
coarser grid, a ratio of $4$ being expected for a second-order accurate method.
The results show that the AMP scheme is stable and close to second-order accurate in the
maximum-norm.

\newcommand{\drawTW}[7]{%
\begin{scope}[#1]
\draw(0,0) node[anchor=south west,xshift=-4pt,yshift=+0pt] {\trimfig{fig/#2}{\figWidth}};
\draw(0,0.0) node[draw,fill=white,anchor=south west,xshift=12pt,yshift=12pt] {\scriptsize #4};
\draw(4.5,.0) node[draw,fill=white,anchor=south east,xshift=-12pt,yshift=12pt] {\scriptsize #5};
\draw (\xcb,\ycb) node[anchor=south west,xshift= +0pt,yshift=+0pt] {\trimfigcb{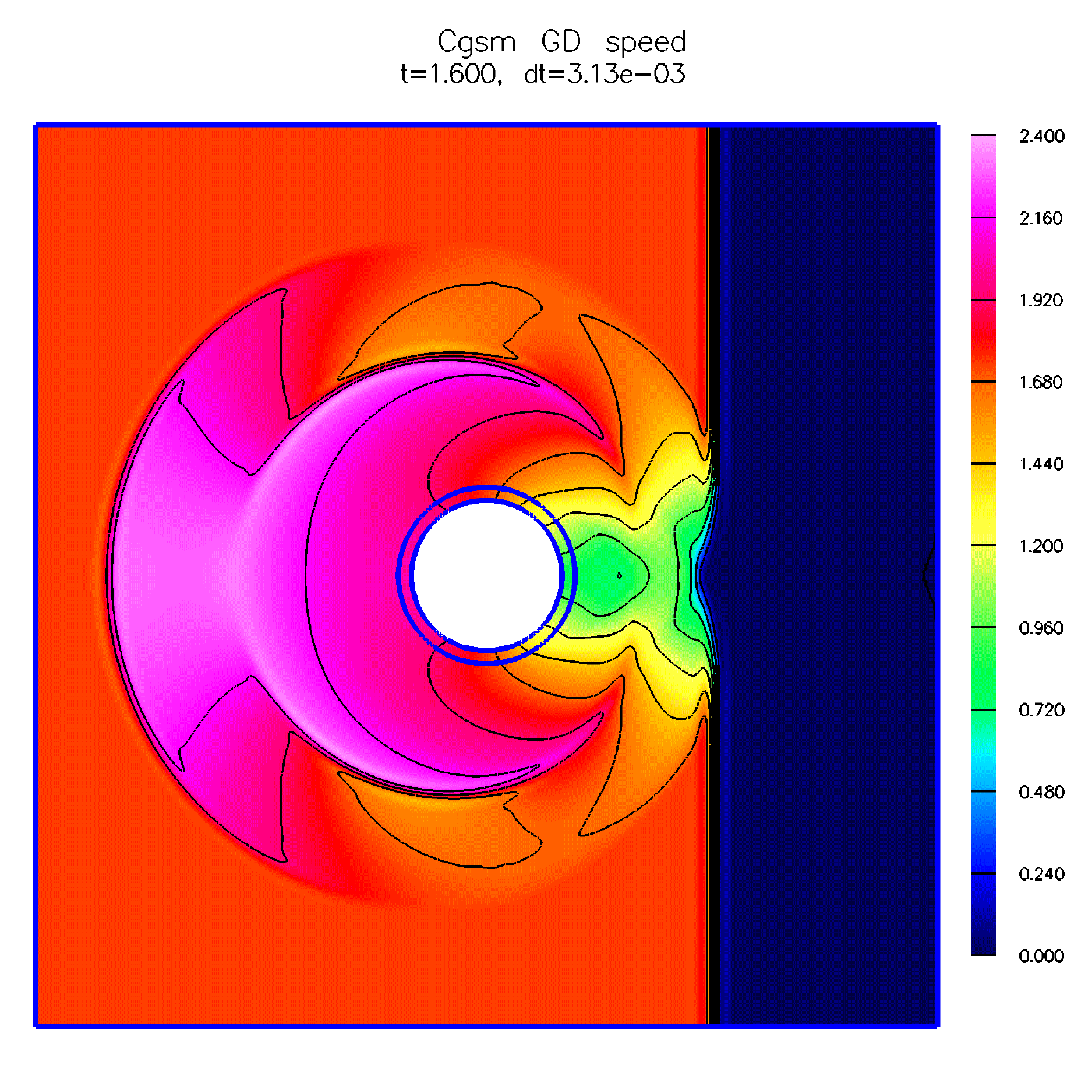}{\cbWidth}{\cbHeight}};
\draw (\xcb,\ycb) node[anchor=south west,xshift= +8pt,yshift=+1pt] {\scriptsize $#6$};
\draw (\xcb,\ycbTop) node[anchor=south west,xshift= +8pt,yshift=-6pt] {\scriptsize $#7$};
\end{scope}
}
{
\newcommand{\cbWidth}{.2cm}
\newcommand{\cbHeight}{4.2cm}
\newcommand{\xcb}{4.4cm}
\newcommand{\ycb}{0.22cm}
\newlength{\ycbTop}%
\setlength{\ycbTop}{\ycb+\cbHeight}
\newlength{\ycbMid}%
\setlength{\ycbMid}{\ycb+\cbHeight*\real{.5}}
\newcommand{\xLabel}{6.5cm}
\newcommand{\yLabel}{6.5cm}
\newcommand{\trimfigcb}[3]{\includegraphics[width=#2, height=#3, clip, trim=17cm 2.35cm 1.65cm 2.35cm]{#1}}
\newcommand{\figWidth}{4.5cm}
\newcommand{\trimfig}[2]{\trimFig{#1}{#2}{.33}{.33}{.25}{.3}}
\begin{figure}[htb]
\begin{center}
\begin{tikzpicture}[scale=1]
  \useasboundingbox (0.2,.75) rectangle (16.4,15.);  
  \drawTW{xshift= -.1cm,yshift=10.0cm}{insShell_mu0_rhos1em2_v1_t1p0}{$v_1$}{$v_1$}{MP-I1, $\delta=10^{-2}$}{-.15}{+.15};
  \drawTW{xshift= 5.5cm,yshift=10.0cm}{insShell_mu0_rhos1em2_v2_t1p0}{$v_2$}{$v_2$}{MP-I1, $\delta=10^{-2}$}{-.15}{+.15};
  \drawTW{xshift=11.0cm,yshift=10.0cm}{insShell_mu0_rhos1em2_p_t1p0}{$p$}{$p$}{MP-I1, $\delta=10^{-2}$}{-.037}{+.037};
%
  \drawTW{xshift= -.1cm,yshift= 5.0cm}{insShell_mu05_rhos1ep0_v1_t1p0}{$v_1$}{$v_1$}{MP-V1, $\delta=1$}{-.21}{+.21};
  \drawTW{xshift= 5.5cm,yshift= 5.0cm}{insShell_mu05_rhos1ep0_v2_t1p0}{$v_2$}{$v_2$}{MP-V1, $\delta=1$}{-.42}{+.42};
  \drawTW{xshift=11.0cm,yshift= 5.0cm}{insShell_mu05_rhos1ep0_p_t1p0}{$p$}{$p$}{MP-V1, $\delta=1$}{-.58}{.58};
%
  \drawTW{xshift= -.1cm,yshift= 0.0cm}{insShell_mu05_rhos1em1_vs2_v1_t1p0}{$v_1$}{$v_1$}{MP-V2, $\delta=.1$}{-.053}{+.053};
  \drawTW{xshift= 5.5cm,yshift= 0.0cm}{insShell_mu05_rhos1em1_vs2_v2_t1p0}{$v_2$}{$v_2$}{MP-V2, $\delta=.1$}{-.11}{+.11};
  \drawTW{xshift=11.0cm,yshift= 0.0cm}{insShell_mu05_rhos1em1_vs2_p_t1p0}{$p$}{$p$}{MP-V2, $\delta=.1$}{-.11}{.11};
%
\end{tikzpicture}
\end{center}
  \caption{Traveling wave solution for a viscous incompressible fluid and a structural shell at $t=1$.
         Top row: model problem MP-I1, inviscid fluid and shell with vertical motion only, $\delta=(\rhos\hs)/(\rho H)=10^{-2}$, $\mu=0$. 
        Middle row: model problem MP-V1, viscous fluid and shell with vertical motion only, $\delta=1$, $\mu=.05$.
        Bottom row: model problem MP-V2, viscous fluid and shell with horizontal and vertical motion, $\delta=.1$, $\mu=.05$.
      }
  \label{fig:insShell}
\end{figure}
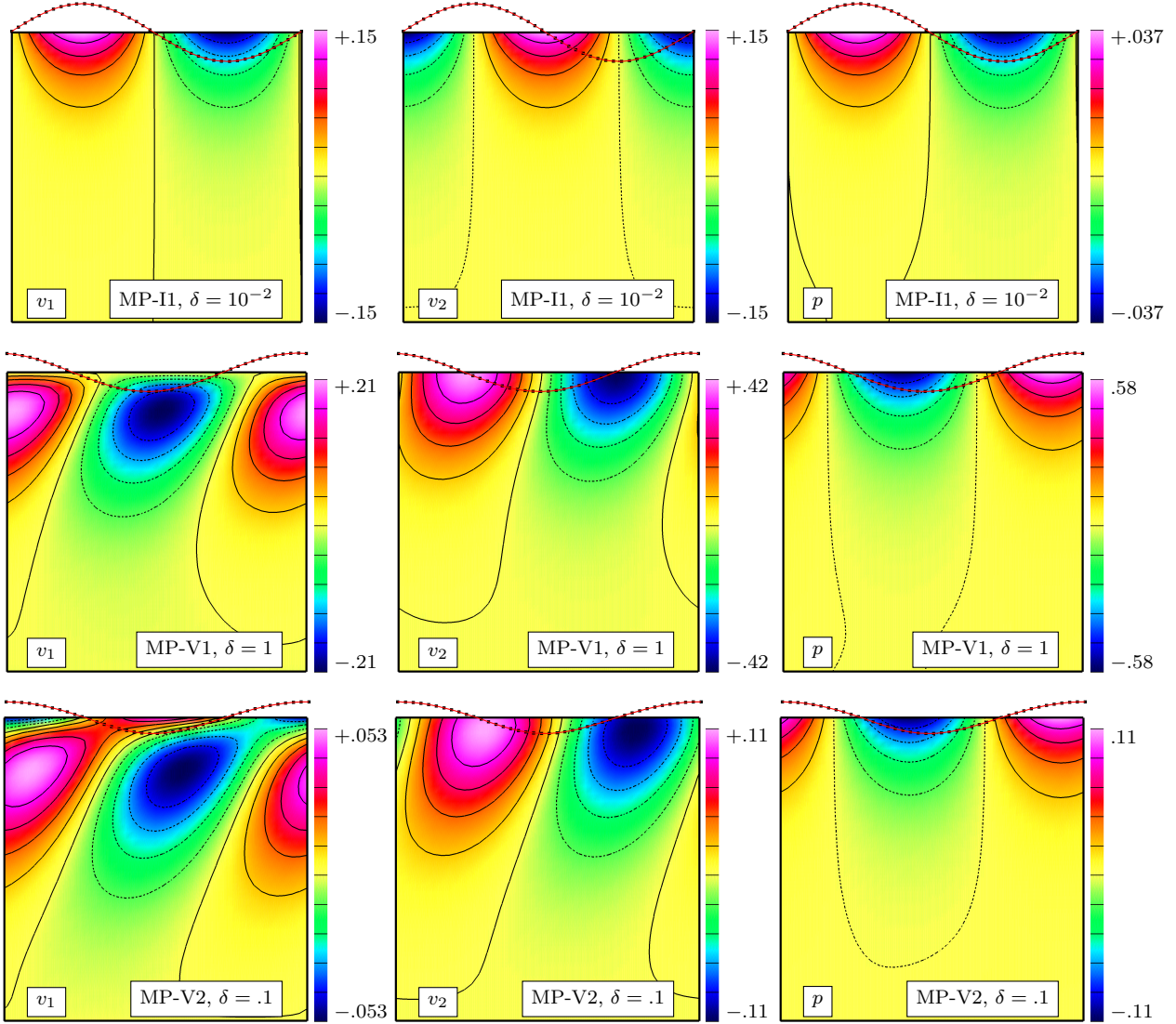
}

\subsection{Traveling wave exact solutions}

Exact traveling wave solutions can be determined for the FSI model problems
defined in Section~\ref{sec:modelProblems}. The traveling wave solutions are of the form
$\vv(x,y,t) = \hat{\vv}(y)e^{i(kx-wt)}$ and $p(x,y,t) = \hat{p}(y)e^{i(kx-wt)}$ for the
velocity and pressure of the fluid, and $\usv(x,t) = \hat{\uv}e^{i(kx-wt)}$ for the displacement of the
shell. The wave is periodic in $x$ with wave number
$k\in\Real$ and has frequency $\omega\in\Complex$ in time.  Solutions for $\hat{\vv}(y)$, $\hat{p}(y)$ and $\hat{\uv}$,
along with a dispersion relation involving $k$ and $\omega$, are given in~\ref{sec:travelingWave} for each
model problem.
  For
viscous fluids the traveling wave decays over time so that ${\rm Im}(\omega)<0$. For
the results shown here, the wave number is taken as $k=2\pi$ and a value of
$\omega$ is chosen so that the wave travels from left to right
(${\rm Re}(\omega)>0$).

The amplitude parameter is chosen as $\ampe=1/10$ which defines the maximum amplitude of the displacement 
on the interface.  The initial conditions are taken from the exact solutions evaluated at $t=0$
and periodic boundary conditions are used in the $x$-direction.

Figure~\ref{fig:insShell} shows shaded contour plots of numerical solutions to the three model problems for a wave traveling from
left to right. 
Contour plots of $v_1$, $v_2$, and $p$ are shown along with the position of the interface for three different
values of the density ratio $\dsf=(\rhos\hs)/(\rho H)$.
The traveling waves are surface waves in the fluid which decay in the vertical direction away from the interface. 
The effects of the fluid viscosity are evident in the components of the velocity, $v_1$ and $v_2$, and in
particular we note the boundary layer in $v_1$ that occurs for the light shell with $\dsf=0.1$.




\begin{figure}[hbt]\tableFontSize
\begin{center}
\begin{tabular}{|c|c|c|c|} \hline
\multicolumn{4}{|c|}{Traveling wave frequencies $\omega$} \\ \hline 
~~$\dsf$~~ & MP-I1        & MP-V1              & MP-V2    \\ \hline 
\strutt $10^{-2}$  &  (1.5277,0)    &  (.25753,$-1.1455$)    &  (.43081,$-1.0018$)          \\ \hline   
\strutt $1$        &  (5.8359,0)       &  (5.6878,$-0.31552$)    &  (5.6467,$-0.34418$)            \\ \hline   
\strutt $10^3$     &  (6.2827,0)       & (6.28253,$-3.8831$e-04) &  (6.2760,$-4.2992$e-03)         \\ \hline   
\end{tabular}
\end{center}
\caption{Values of the (complex) frequency $\omega=a+ib=(a,b)$ for the exact traveling wave solution used in the numerical simulations
of the different model problems for the density ratio $\dsf=(\rhos\hs)/(\rho H)$.}
\label{tab:travelingWaveOmega}
\end{figure}




\bogus{
\begin{table}[hbt]\tableFont 
\begin{center}
\begin{tabular}{|l|c|c|c|c|c|c|c|c|c|} \hline 
grid  & N & $p$ & r &  $\vert \mathbf{v} \vert$ & r & $\bar{\mathbf{u}}$ & r & $\bar{\mathbf{v}}$ & r \\ \hline 
                sq20 &     1 & \num{3.2}{-2} &      & \num{2.3}{-2} &      & \num{2.7}{-3} &      & \num{1.7}{-2} &       \\ \hline
                sq40 &     2 & \num{8.6}{-3} &  3.7 & \num{5.9}{-3} &  3.9 & \num{6.7}{-4} &  4.0 & \num{4.2}{-3} &  4.0  \\ \hline
                sq80 &     4 & \num{2.0}{-3} &  4.3 & \num{1.4}{-3} &  4.2 & \num{1.7}{-4} &  4.0 & \num{1.0}{-3} &  4.0  \\ \hline
               sq160 &     8 & \num{4.3}{-4} &  4.6 & \num{3.2}{-4} &  4.3 & \num{4.2}{-5} &  4.0 & \num{2.6}{-4} &  4.0  \\ \hline
    rate             &       &  $2.07$       &      &  $2.05$       &      &  $2.00$       &      &  $1.99$       &       \\ \hline
\end{tabular}
\caption{Inses shell: fw.mu0.normal.am1.rep3, AM=1, amImp=1, vpf=0.1, ts=pc2, $t=1.$, $\mu=0.0$, $\rho=1$, $\bar{\rho}=1e3$, T=1e3, , cfl=.5, normalMotion=1, cdv=0.25, ad2=1., ad2e=0 gp=pn, Tue Jul 30  7:08:10 2013}\label{table:shell.fw.mu0.normal.am1.rep3.pc2.pc2}
\end{center}
\end{table}
}
\newcommand{\tableMPIIAmpHeavy}{%
\begin{tabular}{|c|c|c|c|c|c|c|c|c|} \hline
\multicolumn{9}{|c|}{MP-I1, traveling wave, heavy solid} \\ \hline 
\strutt~~$h_j$~~& $E_j^{(p)}$ & $r$ & $E_j^{(\vv)}$ & $r$ & $E_j^{(\usv)}$ & $r$  & $E_j^{(\vsv)}$  & $r$  \\ \hline 
 1/20      & \num{3.2}{-2} &      & \num{2.3}{-2} &      & \num{2.7}{-3} &      & \num{1.7}{-2} &       \\ \hline
 1/40      & \num{8.6}{-3} &  3.7 & \num{5.9}{-3} &  3.9 & \num{6.7}{-4} &  4.0 & \num{4.2}{-3} &  4.0  \\ \hline
 1/80      & \num{2.0}{-3} &  4.3 & \num{1.4}{-3} &  4.2 & \num{1.7}{-4} &  4.0 & \num{1.0}{-3} &  4.0  \\ \hline
 1/160	   & \num{4.3}{-4} &  4.6 & \num{3.2}{-4} &  4.3 & \num{4.2}{-5} &  4.0 & \num{2.6}{-4} &  4.0  \\ \hline
\rateLabel &  $2.07$       &      &  $2.05$       &      &  $2.00$       &      &  $1.99$       &       \\ \hline
\end{tabular}
}

\bogus{
\begin{table}[hbt]\tableFont 
\begin{center}
\begin{tabular}{|l|c|c|c|c|c|c|c|c|c|} \hline 
grid  & N & $p$ & r &  $\vert \mathbf{v} \vert$ & r & $\bar{\mathbf{u}}$ & r & $\bar{\mathbf{v}}$ & r \\ \hline 
                sq20 &     1 & \num{2.6}{-2} &      & \num{2.7}{-2} &      & \num{3.8}{-3} &      & \num{2.0}{-2} &       \\ \hline
                sq40 &     2 & \num{5.8}{-3} &  4.5 & \num{5.5}{-3} &  5.0 & \num{8.6}{-4} &  4.5 & \num{4.6}{-3} &  4.4  \\ \hline
                sq80 &     4 & \num{1.3}{-3} &  4.4 & \num{1.1}{-3} &  4.7 & \num{2.1}{-4} &  4.1 & \num{1.1}{-3} &  4.1  \\ \hline
               sq160 &     8 & \num{3.0}{-4} &  4.3 & \num{2.8}{-4} &  4.2 & \num{5.1}{-5} &  4.1 & \num{2.8}{-4} &  4.0  \\ \hline
    rate             &       &  $2.14$       &      &  $2.21$       &      &  $2.07$       &      &  $2.07$       &       \\ \hline
\end{tabular}
\caption{Inses shell: fw.mu0.normal.am1.rep0, AM=1, amImp=1, vpf=0.1, ts=pc2, $t=1.$, $\mu=0.0$, $\rho=1$, $\bar{\rho}=1$, T=1, , cfl=.5, normalMotion=1, cdv=0.25, ad2=2., ad2e=0 gp=pn, Tue Jul 30  7:12:14 2013}\label{table:shell.fw.mu0.normal.am1.rep0.pc2.pc2}
\end{center}
\end{table}
}
\newcommand{\tableMPIIAmpMedium}{%
\begin{tabular}{|c|c|c|c|c|c|c|c|c|} \hline
\multicolumn{9}{|c|}{MP-I1, traveling wave, medium solid} \\ \hline 
\strutt~~$h_j$~~& $E_j^{(p)}$ & $r$ & $E_j^{(\vv)}$ & $r$ & $E_j^{(\usv)}$ & $r$  & $E_j^{(\vsv)}$  & $r$  \\ \hline 
 1/20      & \num{2.6}{-2} &      & \num{2.7}{-2} &      & \num{3.8}{-3} &      & \num{2.0}{-2} &       \\ \hline
 1/40      & \num{5.8}{-3} &  4.5 & \num{5.5}{-3} &  5.0 & \num{8.6}{-4} &  4.5 & \num{4.6}{-3} &  4.4  \\ \hline
 1/80      & \num{1.3}{-3} &  4.4 & \num{1.1}{-3} &  4.7 & \num{2.1}{-4} &  4.1 & \num{1.1}{-3} &  4.1  \\ \hline
 1/160	   & \num{3.0}{-4} &  4.3 & \num{2.8}{-4} &  4.2 & \num{5.1}{-5} &  4.1 & \num{2.8}{-4} &  4.0  \\ \hline
\rateLabel &  $2.14$       &      &  $2.21$       &      &  $2.07$       &      &  $2.07$       &       \\ \hline
\end{tabular}
}

\bogus{
\begin{table}[hbt]\tableFont 
\begin{center}
\begin{tabular}{|l|c|c|c|c|c|c|c|c|c|} \hline 
grid  & N & $p$ & r &  $\vert \mathbf{v} \vert$ & r & $\bar{\mathbf{u}}$ & r & $\bar{\mathbf{v}}$ & r \\ \hline 
                sq20 &     1 & \num{7.7}{-4} &      & \num{4.4}{-3} &      & \num{1.6}{-3} &      & \num{2.8}{-3} &       \\ \hline
                sq40 &     2 & \num{1.9}{-4} &  4.1 & \num{1.2}{-3} &  3.6 & \num{3.9}{-4} &  4.0 & \num{7.8}{-4} &  3.6  \\ \hline
                sq80 &     4 & \num{4.5}{-5} &  4.1 & \num{3.4}{-4} &  3.7 & \num{9.9}{-5} &  4.0 & \num{2.1}{-4} &  3.7  \\ \hline
               sq160 &     8 & \num{1.1}{-5} &  4.2 & \num{8.8}{-5} &  3.8 & \num{2.3}{-5} &  4.2 & \num{5.6}{-5} &  3.8  \\ \hline
    rate             &       &  $2.05$       &      &  $1.89$       &      &  $2.02$       &      &  $1.88$       &       \\ \hline
\end{tabular}
\caption{Inses shell: fw.mu0.normal.am1.rem2, AM=1, amImp=1, vpf=0.1, ts=pc2, $t=1.$, $\mu=0.0$, $\rho=1$, $\bar{\rho}=1e-2$, T=1e-2, , cfl=.5, normalMotion=1, cdv=0.25, ad2=1., ad2e=0 gp=pn, Tue Jul 30  6:59:00 2013}\label{table:shell.fw.mu0.normal.am1.rem2.pc2.pc2}
\end{center}
\end{table}
}
\newcommand{\tableMPIIAmpLight}{%
\begin{tabular}{|c|c|c|c|c|c|c|c|c|} \hline
\multicolumn{9}{|c|}{MP-I1, traveling wave, light solid} \\ \hline 
\strutt~~$h_j$~~& $E_j^{(p)}$ & $r$ & $E_j^{(\vv)}$ & $r$ & $E_j^{(\usv)}$ & $r$  & $E_j^{(\vsv)}$  & $r$  \\ \hline 
 1/20      & \num{7.7}{-4} &      & \num{4.4}{-3} &      & \num{1.6}{-3} &      & \num{2.8}{-3} &       \\ \hline 
 1/40      & \num{1.9}{-4} &  4.1 & \num{1.2}{-3} &  3.6 & \num{3.9}{-4} &  4.0 & \num{7.8}{-4} &  3.6  \\ \hline 
 1/80      & \num{4.5}{-5} &  4.1 & \num{3.4}{-4} &  3.7 & \num{9.9}{-5} &  4.0 & \num{2.1}{-4} &  3.7  \\ \hline 
 1/160	   & \num{1.1}{-5} &  4.2 & \num{8.8}{-5} &  3.8 & \num{2.3}{-5} &  4.2 & \num{5.6}{-5} &  3.8  \\ \hline 
\rateLabel &  $2.05$       &      &  $1.89$       &      &  $2.02$       &      &  $1.88$       &       \\ \hline 
\end{tabular}
}

\begin{figure}[hbt]\tableFontSize
\begin{center}
\tableMPIIAmpHeavy
\vskip.5\baselineskip
\tableMPIIAmpMedium
\vskip.5\baselineskip
\tableMPIIAmpLight
\caption{Traveling wave solution for an inviscid incompressible fluid and a shell that only supports
vertical motion (model problem MP-I1).
Maximum errors and estimated convergence rates at $t=1.0$, computed using the AMP scheme for a heavy solid, $\delta=10^3$,
medium solid, $\delta=1$, and light solid, $\delta=10^{-2}$, where $\delta=(\rhos\hs)/(\rho H)$.
}
\label{tab:TWMaxNormInviscid}
\end{center}
\end{figure}

\bogus{
\begin{table}[hbt]\tableFont 
\begin{center}
\begin{tabular}{|l|c|c|c|c|c|c|c|c|c|} \hline 
grid  & N & $p$ & r &  $\vert \mathbf{v} \vert$ & r & $\bar{\mathbf{u}}$ & r & $\bar{\mathbf{v}}$ & r \\ \hline 
                sq20 &     1 & \num{8.1}{-2} &      & \num{2.3}{-2} &      & \num{1.7}{-2} &      & \num{2.3}{-2} &       \\ \hline
                sq40 &     2 & \num{1.2}{-2} &  6.5 & \num{5.1}{-3} &  4.5 & \num{4.3}{-3} &  3.9 & \num{5.1}{-3} &  4.5  \\ \hline
                sq80 &     4 & \num{2.5}{-3} &  5.0 & \num{1.1}{-3} &  4.5 & \num{1.0}{-3} &  4.1 & \num{1.1}{-3} &  4.5  \\ \hline
               sq160 &     8 & \num{6.0}{-4} &  4.2 & \num{2.7}{-4} &  4.2 & \num{2.6}{-4} &  4.1 & \num{2.7}{-4} &  4.2  \\ \hline
    rate             &       &  $2.36$       &      &  $2.14$       &      &  $2.01$       &      &  $2.14$       &       \\ \hline
\end{tabular}
\caption{Inses shell: fw.mup05.normal.am1.rep3, AM=1, amImp=1, vpf=0.1, ts=pc2, $t=.5$, $\mu=0.05$, $\rho=1$, $\bar{\rho}=1e3$, T=1e3, , cfl=0.75, normalMotion=1, ad2=0, ad2e=0 gp=pn, Mon Jul 29 11:24:46 2013}\label{table:shell.fw.mup05.normal.am1.rep3.pc2.pc2}
\end{center}
\end{table}
}
\newcommand{\tableMPVIAmpHeavy}{%
\begin{tabular}{|c|c|c|c|c|c|c|c|c|} \hline
\multicolumn{9}{|c|}{MP-V1, traveling wave, heavy solid} \\ \hline 
\strutt~~$h_j$~~& $E_j^{(p)}$ & $r$ & $E_j^{(\vv)}$ & $r$ & $E_j^{(\usv)}$ & $r$  & $E_j^{(\vsv)}$  & $r$  \\ \hline 
 1/20      &  \num{8.1}{-2} &      & \num{2.3}{-2} &      & \num{1.7}{-2} &      & \num{2.3}{-2} &       \\ \hline
 1/40      &  \num{1.2}{-2} &  6.5 & \num{5.1}{-3} &  4.5 & \num{4.3}{-3} &  3.9 & \num{5.1}{-3} &  4.5  \\ \hline
 1/80      &  \num{2.5}{-3} &  5.0 & \num{1.1}{-3} &  4.5 & \num{1.0}{-3} &  4.1 & \num{1.1}{-3} &  4.5  \\ \hline
 1/160	   &  \num{6.0}{-4} &  4.2 & \num{2.7}{-4} &  4.2 & \num{2.6}{-4} &  4.1 & \num{2.7}{-4} &  4.2  \\ \hline
\rateLabel &   $2.36$       &      &  $2.14$       &      &  $2.01$       &      &  $2.14$       &       \\ \hline
\end{tabular}
}

\bogus{
\begin{table}[hbt]\tableFont 
\begin{center}
\begin{tabular}{|l|c|c|c|c|c|c|c|c|c|} \hline 
grid  & N & $p$ & r &  $\vert \mathbf{v} \vert$ & r & $\bar{\mathbf{u}}$ & r & $\bar{\mathbf{v}}$ & r \\ \hline 
                sq20 &     1 & \num{1.6}{-2} &      & \num{1.2}{-2} &      & \num{7.1}{-3} &      & \num{1.2}{-2} &       \\ \hline
                sq40 &     2 & \num{3.5}{-3} &  4.4 & \num{2.6}{-3} &  4.8 & \num{1.8}{-3} &  4.0 & \num{2.6}{-3} &  4.8  \\ \hline
                sq80 &     4 & \num{7.9}{-4} &  4.4 & \num{5.6}{-4} &  4.6 & \num{4.2}{-4} &  4.3 & \num{5.6}{-4} &  4.6  \\ \hline
               sq160 &     8 & \num{1.9}{-4} &  4.2 & \num{1.3}{-4} &  4.2 & \num{1.0}{-4} &  4.2 & \num{1.3}{-4} &  4.2  \\ \hline
    rate             &       &  $2.12$       &      &  $2.19$       &      &  $2.06$       &      &  $2.19$       &       \\ \hline
\end{tabular}
\caption{Inses shell: fw.mup05.normal.am1.rep0, AM=1, amImp=1, vpf=0.1, ts=pc2, $t=.5$, $\mu=0.05$, $\rho=1$, $\bar{\rho}=1$, T=1, , cfl=0.75, normalMotion=1, ad2=0, ad2e=0 gp=pn, Mon Jul 29 12:31:31 2013}\label{table:shell.fw.mup05.normal.am1.rep0.pc2.pc2}
\end{center}
\end{table}
}
\newcommand{\tableMPVIAmpMedium}{%
\begin{tabular}{|c|c|c|c|c|c|c|c|c|} \hline
\multicolumn{9}{|c|}{MP-V1, traveling wave, medium solid} \\ \hline 
\strutt~~$h_j$~~& $E_j^{(p)}$ & $r$ & $E_j^{(\vv)}$ & $r$ & $E_j^{(\usv)}$ & $r$  & $E_j^{(\vsv)}$  & $r$  \\ \hline 
 1/20      &  \num{1.6}{-2} &      & \num{1.2}{-2} &      & \num{7.1}{-3} &      & \num{1.2}{-2} &       \\ \hline
 1/40      &  \num{3.5}{-3} &  4.4 & \num{2.6}{-3} &  4.8 & \num{1.8}{-3} &  4.0 & \num{2.6}{-3} &  4.8  \\ \hline
 1/80      &  \num{7.9}{-4} &  4.4 & \num{5.6}{-4} &  4.6 & \num{4.2}{-4} &  4.3 & \num{5.6}{-4} &  4.6  \\ \hline
 1/160	   &  \num{1.9}{-4} &  4.2 & \num{1.3}{-4} &  4.2 & \num{1.0}{-4} &  4.2 & \num{1.3}{-4} &  4.2  \\ \hline
\rateLabel &   $2.12$       &      &  $2.19$       &      &  $2.06$       &      &  $2.19$       &       \\ \hline
\end{tabular}
}

\bogus{
\begin{table}[hbt]\tableFont 
\begin{center}
\begin{tabular}{|l|c|c|c|c|c|c|c|c|c|} \hline 
grid  & N & $p$ & r &  $\vert \mathbf{v} \vert$ & r & $\bar{\mathbf{u}}$ & r & $\bar{\mathbf{v}}$ & r \\ \hline 
                sq20 &     1 & \num{3.4}{-4} &      & \num{2.2}{-4} &      & \num{9.3}{-4} &      & \num{9.3}{-4} &       \\ \hline
                sq40 &     2 & \num{8.2}{-5} &  4.1 & \num{7.4}{-5} &  3.0 & \num{2.1}{-4} &  4.5 & \num{2.1}{-4} &  4.5  \\ \hline
                sq80 &     4 & \num{2.0}{-5} &  4.1 & \num{2.0}{-5} &  3.7 & \num{4.9}{-5} &  4.3 & \num{4.9}{-5} &  4.3  \\ \hline
               sq160 &     8 & \num{5.0}{-6} &  4.0 & \num{5.1}{-6} &  3.9 & \num{1.2}{-5} &  4.1 & \num{1.2}{-5} &  4.1  \\ \hline
    rate             &       &  $2.03$       &      &  $1.81$       &      &  $2.10$       &      &  $2.10$       &       \\ \hline
\end{tabular}
\caption{Inses shell: fw.mup05.normal.am1.rem2, AM=1, amImp=1, vpf=0.1, ts=pc2, $t=.5$, $\mu=0.05$, $\rho=1$, $\bar{\rho}=1e-2$, T=1e-2, , cfl=0.75, normalMotion=1, ad2=0, ad2e=0 gp=pn, Mon Jul 29 15:02:45 2013}\label{table:shell.fw.mup05.normal.am1.rem2.pc2.pc2}
\end{center}
\end{table}
}
\newcommand{\tableMPVIAmpLight}{%
\begin{tabular}{|c|c|c|c|c|c|c|c|c|} \hline
\multicolumn{9}{|c|}{MP-V1, traveling wave, light solid} \\ \hline 
\strutt~~$h_j$~~& $E_j^{(p)}$ & $r$ & $E_j^{(\vv)}$ & $r$ & $E_j^{(\usv)}$ & $r$  & $E_j^{(\vsv)}$  & $r$  \\ \hline 
 1/20      &  \num{3.4}{-4} &      & \num{2.2}{-4} &      & \num{9.3}{-4} &      & \num{9.3}{-4} &       \\ \hline
 1/40      &  \num{8.2}{-5} &  4.1 & \num{7.4}{-5} &  3.0 & \num{2.1}{-4} &  4.5 & \num{2.1}{-4} &  4.5  \\ \hline
 1/80      &  \num{2.0}{-5} &  4.1 & \num{2.0}{-5} &  3.7 & \num{4.9}{-5} &  4.3 & \num{4.9}{-5} &  4.3  \\ \hline
 1/160	   &  \num{5.0}{-6} &  4.0 & \num{5.1}{-6} &  3.9 & \num{1.2}{-5} &  4.1 & \num{1.2}{-5} &  4.1  \\ \hline
\rateLabel &   $2.03$       &      &  $1.81$       &      &  $2.10$       &      &  $2.10$       &       \\ \hline
\end{tabular}
}

\begin{figure}[hbt]\tableFontSize
\begin{center}
\tableMPVIAmpHeavy
\vskip.5\baselineskip
\tableMPVIAmpMedium
\vskip.5\baselineskip
\tableMPVIAmpLight
\caption{Traveling wave solution for a viscous incompressible fluid and a shell that only supports
vertical motion (model problem MP-V1).
Maximum errors and estimated convergence rates at $t=0.5$, computed using the AMP scheme for a heavy solid, $\delta=10^3$,
medium solid, $\delta=1$, and light solid, $\delta=10^{-2}$, where $\delta=(\rhos\hs)/(\rho H)$.
}
\label{tab:TWMaxNormViscousNormalMotion}
\end{center}
\end{figure}

\bogus{
\begin{table}[hbt]\tableFont 
\begin{center}
\begin{tabular}{|l|c|c|c|c|c|c|c|c|c|} \hline 
grid  & N & $p$ & r &  $\vert \mathbf{v} \vert$ & r & $\bar{\mathbf{u}}$ & r & $\bar{\mathbf{v}}$ & r \\ \hline 
                sq20 &     1 & \num{6.0}{-2} &      & \num{2.7}{-2} &      & \num{2.0}{-2} &      & \num{2.7}{-2} &       \\ \hline
                sq40 &     2 & \num{1.3}{-2} &  4.5 & \num{6.0}{-3} &  4.4 & \num{5.0}{-3} &  4.0 & \num{6.0}{-3} &  4.4  \\ \hline
                sq80 &     4 & \num{3.0}{-3} &  4.4 & \num{1.4}{-3} &  4.4 & \num{1.2}{-3} &  4.1 & \num{1.4}{-3} &  4.4  \\ \hline
               sq160 &     8 & \num{7.3}{-4} &  4.1 & \num{3.2}{-4} &  4.2 & \num{3.0}{-4} &  4.1 & \num{3.2}{-4} &  4.2  \\ \hline
    rate             &       &  $2.12$       &      &  $2.12$       &      &  $2.02$       &      &  $2.12$       &       \\ \hline
\end{tabular}
\caption{Inses shell: fw.mup05.full.am1.rep2, AM=1, amImp=1, vpf=0.1, ts=pc2, $t=.5$, $\mu=0.05$, $\rho=1$, $\bar{\rho}=1e2$, T=1e2, , cfl=0.75, normalMotion=0, gp=pn, Mon Jul 29  8:14:15 2013}\label{table:shell.fw.mup05.full.am1.rep2.pc2.pc2}
\end{center}
\end{table}
}
\newcommand{\tableMPVIIAmpHeavy}{%
\begin{tabular}{|c|c|c|c|c|c|c|c|c|} \hline
\multicolumn{9}{|c|}{MP-V2, traveling wave, heavy solid} \\ \hline 
\strutt~~$h_j$~~& $E_j^{(p)}$ & $r$ & $E_j^{(\vv)}$ & $r$ & $E_j^{(\usv)}$ & $r$  & $E_j^{(\vsv)}$  & $r$  \\ \hline 
 1/20      &  \num{6.0}{-2} &      & \num{2.7}{-2} &      & \num{2.0}{-2} &      & \num{2.7}{-2} &       \\ \hline
 1/40      &  \num{1.3}{-2} &  4.5 & \num{6.0}{-3} &  4.4 & \num{5.0}{-3} &  4.0 & \num{6.0}{-3} &  4.4  \\ \hline
 1/80      &  \num{3.0}{-3} &  4.4 & \num{1.4}{-3} &  4.4 & \num{1.2}{-3} &  4.1 & \num{1.4}{-3} &  4.4  \\ \hline
 1/160	   &  \num{7.3}{-4} &  4.1 & \num{3.2}{-4} &  4.2 & \num{3.0}{-4} &  4.1 & \num{3.2}{-4} &  4.2  \\ \hline
\rateLabel &   $2.12$       &      &  $2.12$       &      &  $2.02$       &      &  $2.12$       &       \\ \hline
\end{tabular}
}

\bogus{
\begin{table}[hbt]\tableFont 
\begin{center}
\begin{tabular}{|l|c|c|c|c|c|c|c|c|c|} \hline 
grid  & N & $p$ & r &  $\vert \mathbf{v} \vert$ & r & $\bar{\mathbf{u}}$ & r & $\bar{\mathbf{v}}$ & r \\ \hline 
                sq20 &     1 & \num{2.0}{-2} &      & \num{1.5}{-2} &      & \num{9.8}{-3} &      & \num{1.5}{-2} &       \\ \hline
                sq40 &     2 & \num{4.6}{-3} &  4.2 & \num{3.2}{-3} &  4.7 & \num{2.4}{-3} &  4.1 & \num{3.2}{-3} &  4.7  \\ \hline
                sq80 &     4 & \num{1.1}{-3} &  4.2 & \num{7.0}{-4} &  4.5 & \num{5.7}{-4} &  4.2 & \num{7.0}{-4} &  4.5  \\ \hline
               sq160 &     8 & \num{2.7}{-4} &  4.1 & \num{1.7}{-4} &  4.2 & \num{1.4}{-4} &  4.1 & \num{1.7}{-4} &  4.2  \\ \hline
    rate             &       &  $2.06$       &      &  $2.17$       &      &  $2.05$       &      &  $2.17$       &       \\ \hline
\end{tabular}
\caption{Inses shell: fw.mup05.full.am1.rep0, AM=1, amImp=1, vpf=0.1, ts=pc2, $t=.5$, $\mu=0.05$, $\rho=1$, $\bar{\rho}=1.$, T=1., , cfl=0.75, normalMotion=0, gp=pn, Mon Jul 29  8:21:12 2013}\label{table:shell.fw.mup05.full.am1.rep0.pc2.pc2}
\end{center}
\end{table}
}
\newcommand{\tableMPVIIAmpMedium}{%
\begin{tabular}{|c|c|c|c|c|c|c|c|c|} \hline
\multicolumn{9}{|c|}{MP-V2, traveling wave, medium solid} \\ \hline 
\strutt~~$h_j$~~& $E_j^{(p)}$ & $r$ & $E_j^{(\vv)}$ & $r$ & $E_j^{(\usv)}$ & $r$  & $E_j^{(\vsv)}$  & $r$  \\ \hline 
 1/20      &  \num{2.0}{-2} &      & \num{1.5}{-2} &      & \num{9.8}{-3} &      & \num{1.5}{-2} &       \\ \hline
 1/40      &  \num{4.6}{-3} &  4.2 & \num{3.2}{-3} &  4.7 & \num{2.4}{-3} &  4.1 & \num{3.2}{-3} &  4.7  \\ \hline
 1/80      &  \num{1.1}{-3} &  4.2 & \num{7.0}{-4} &  4.5 & \num{5.7}{-4} &  4.2 & \num{7.0}{-4} &  4.5  \\ \hline
 1/160	   &  \num{2.7}{-4} &  4.1 & \num{1.7}{-4} &  4.2 & \num{1.4}{-4} &  4.1 & \num{1.7}{-4} &  4.2  \\ \hline
\rateLabel &   $2.06$       &      &  $2.17$       &      &  $2.05$       &      &  $2.17$       &       \\ \hline
\end{tabular}
}

\bogus{
\begin{table}[hbt]\tableFont 
\begin{center}
\begin{tabular}{|l|c|c|c|c|c|c|c|c|c|} \hline 
grid  & N & $p$ & r &  $\vert \mathbf{v} \vert$ & r & $\bar{\mathbf{u}}$ & r & $\bar{\mathbf{v}}$ & r \\ \hline 
                sq20 &     1 & \num{7.4}{-4} &      & \num{5.4}{-4} &      & \num{1.6}{-3} &      & \num{1.6}{-3} &       \\ \hline
                sq40 &     2 & \num{2.1}{-4} &  3.4 & \num{9.6}{-5} &  5.6 & \num{3.8}{-4} &  4.4 & \num{3.8}{-4} &  4.4  \\ \hline
                sq80 &     4 & \num{5.7}{-5} &  3.8 & \num{2.0}{-5} &  4.9 & \num{8.9}{-5} &  4.2 & \num{8.9}{-5} &  4.2  \\ \hline
               sq160 &     8 & \num{1.4}{-5} &  3.9 & \num{5.1}{-6} &  3.9 & \num{2.2}{-5} &  4.1 & \num{2.2}{-5} &  4.1  \\ \hline
    rate             &       &  $1.90$       &      &  $2.25$       &      &  $2.08$       &      &  $2.08$       &       \\ \hline
\end{tabular}
\caption{Inses shell: fw.mup05.full.am1.rem2, AM=1, amImp=1, vpf=0.1, ts=pc2, $t=.5$, $\mu=0.05$, $\rho=1$, $\bar{\rho}=1e-2$, T=1e-2, , cfl=0.75, normalMotion=0, gp=pn, Mon Jul 29  8:35:01 2013}\label{table:shell.fw.mup05.full.am1.rem2.pc2.pc2}
\end{center}
\end{table}
}
\newcommand{\tableMPVIIAmpLight}{%
\begin{tabular}{|c|c|c|c|c|c|c|c|c|} \hline
\multicolumn{9}{|c|}{MP-V2, traveling wave, light solid} \\ \hline 
\strutt~~$h_j$~~& $E_j^{(p)}$ & $r$ & $E_j^{(\vv)}$ & $r$ & $E_j^{(\usv)}$ & $r$  & $E_j^{(\vsv)}$  & $r$  \\ \hline 
 1/20      &  \num{7.4}{-4} &      & \num{5.4}{-4} &      & \num{1.6}{-3} &      & \num{1.6}{-3} &       \\ \hline
 1/40      &  \num{2.1}{-4} &  3.4 & \num{9.6}{-5} &  5.6 & \num{3.8}{-4} &  4.4 & \num{3.8}{-4} &  4.4  \\ \hline
 1/80      &  \num{5.7}{-5} &  3.8 & \num{2.0}{-5} &  4.9 & \num{8.9}{-5} &  4.2 & \num{8.9}{-5} &  4.2  \\ \hline
 1/160	   &  \num{1.4}{-5} &  3.9 & \num{5.1}{-6} &  3.9 & \num{2.2}{-5} &  4.1 & \num{2.2}{-5} &  4.1  \\ \hline
\rateLabel &   $1.90$       &      &  $2.25$       &      &  $2.08$       &      &  $2.08$       &       \\ \hline
\end{tabular}
}

\begin{figure}[hbt]\tableFontSize
\begin{center}
\tableMPVIIAmpHeavy
\vskip.5\baselineskip
\tableMPVIIAmpMedium
\vskip.5\baselineskip
\tableMPVIIAmpLight
\caption{Traveling wave solution for a viscous incompressible fluid and structural shell (model problem MP-V2).
Maximum errors and estimated convergence rates at $t=0.5$, computed using the AMP scheme for a heavy solid, $\delta=10^3$,
medium solid, $\delta=1$, and light solid, $\delta=10^{-2}$, where $\delta=(\rhos\hs)/(\rho H)$.
}
\label{tab:TWMaxNormViscousFullMotion}
\end{center}
\end{figure}

The accuracy and stability of the AMP algorithm for the different model problems can be evaluated quantitatively by
computing traveling wave solutions, for different density ratios, on a sequence of grids of increasing resolution.
The frequencies $\omega$ for the corresponding exact solutions are given in Figure~\ref{tab:travelingWaveOmega}.
The results of this study are summarized in the tables found in Figures~\ref{tab:TWMaxNormInviscid}--\ref{tab:TWMaxNormViscousFullMotion}.
For each model problem, MP-I1, MP-V1, and MP-V2, max-norm errors and estimated convergence rates
are given for density ratios $\dsf=10^{-2}$, $1$ and $10^3$.
The results show that the AMP predictor-corrector scheme is stable and close to second-order accurate in the max-norm. 
Recall that for model problem MP-I1, the fluid viscosity is zero. 
In this case, since the discretization of the fluid equations uses central finite differences, a small
amount of artificial dissipation, proportional to $h_j^2$, is added to the fluid momentum equation,
in order to smooth boundary layers in the error that otherwise degrade the max-norm convergence
rates somewhat (the scheme is stable without dissipation).

\section{Conclusions} \label{sec:conclusions}

We have described a stable added-mass partitioned (AMP) algorithm for the
solution of FSI problems that couple incompressible flows to structural shells (or beams).
The scheme was developed and evaluated for a linearized problem where the fluid
is modeled with the Stokes equations on a fixed reference domain
and the structure is modeled with a linear beam or generalized string model.
The AMP algorithm is based on generalized Robin boundary conditions for the fluid
pressure and velocity that  were derived at a continuous level 
by combining the fluid momentum equation with the shell equation. The conditions
can also be derived at a fully discrete level which ensures that the fluid and structure velocities match exactly on
the interface; otherwise a density weighted projection is used to enforce matching of the velocities.
In addition, suitable forms of the AMP conditions
for fractional-step fluid solvers were provided.
Using mode analysis, the stability of the AMP scheme was proved for a two-dimensional model problem. 
The analysis showed that the AMP scheme is stable even for very light structures and requires no sub-iterations. 
A second-order accurate predictor-corrector algorithm that implements the AMP scheme was described. 
Numerical results for three model FSI model problems in two dimensions were obtained and compared with
exact solutions constructed using the method of analytic solutions and exact traveling wave solutions.
These results verified the stability of the AMP algorithm and demonstated the second-order accuracy of the
scheme in the maxiumum-norm.

In future work, the AMP conditions will be incorporated into fully nonlinear FSI schemes that
treat large displacements and deformations and are based on the deforming composite grid approach~\cite{fsi2012}.
The scheme will also be extended to more general beam and plate models for the structure.

\appendix

\renewcommand{\usvh}{\uvh}
\section{Traveling wave exact solutions for the  FSI model problems} \label{sec:travelingWave}

Exact {\em traveling wave} solutions of the three model
problems, MP-I1, MP-V1 and MP-V2, defined in Section~\ref{sec:modelProblems} are described briefly.
The traveling wave solutions have the general form
\[
\vv(x,y,t)=\vvh(y)\expkw, \qquad p(x,y,t)=\ph(y)\expkw, \qquad \usv(x,t) = \usvh\,\expkw, 
\]
where $k$ is a wave number, $\omega$ is a frequency (possibly complex) and $\vvh(y)$, $\ph(y)$ and $\usvh$ are to be determined.  The solutions are assumed to be $2\pi$-periodic in the $x$ direction so that $k$ takes integer values.

\subsection{Traveling wave solution for MP-I1}

The model problem MP-I1 involves an inviscid incompressible fluid and a shell that supports motion in the vertical direction only.  An exact traveling wave solution is given by 
\begin{equation}
\begin{array}{l}
\displaystyle{
 \eta(x,t) =  \ampe \,\expkw
} \smallskip\\
\displaystyle{
 v_1(x,y,t) = \ampe  \frac{\omega\cosh(k (y+H))}{\sinh(k H)} \,\expkw
} \smallskip\\
\displaystyle{
 v_2(x,y,t) = -\ampe \frac{i \omega\sinh(k (y+H))}{\sinh(k H)} \,\expkw
} \smallskip\\
\displaystyle{
   p(x,y,t) = \ampe \frac{\rho \omega^2\cosh(k (y+H))}{k\sinh(k H)} \,\expkw
}
\end{array}
\label{eq:twsoln}
\end{equation}
where $\eta(x,t)=\bar u_2(x,t)$ is the vertical displacement of the shell, and $\omega$ satisfies the dispersion relation
\begin{align}
  \omega = W(k) = \pm \sqrt{ \frac{\Ks + k^2 \Ts}{ \rhos \hs+ M_a} }. \label{eq:dispersionMP-I1} 
\end{align}
Recall that $\Ks$ is a stiffness coefficient, $\Ts$ is a coefficient of tension, $\rhos$ is the density and $\hs$ is the thickness of the shell.  Here, $M_a$ is the coefficient of added-mass given by 
\begin{align}
&  M_a = \frac{\rho}{k\tanh(k H)},   \label{eq:inviscidAM}
\end{align}
where $\rho$ is the density of the fluid and $H$ is the depth of the fluid domain.
The amplitude parameter $\ampe$ is an arbitrary constant. The real and imaginary parts of~\eqref{eq:twsoln} define real solutions to the model problem MP-I1. 

The dispersion relation~\eqref{eq:dispersionMP-I1} shows that the effect of the fluid
on the shell, as given by $M_a>0$, reduces $\omega$ and the phase speed, $\omega/k$, of the wave. 
Furthermore, \eqref{eq:inviscidAM} indicates that the effect of the added mass increase when $\rho$ increases,
or when $k$ or $H$ decrease.

\vskip\baselineskip


\subsection{Traveling wave solutions for MP-V1 and MP-V2}

The model problems MP-V1 and MP-V2 involve a viscous incompressible fluid, with viscosity $\mu$, coupled to a shell that supports either motion in the vertical direction only (MP-V1) or motion in both the horizontal and vertical directions (MP-V2).  Traveling wave solutions satisfying the governing equations in the fluid domain have the form
\begin{equation}
\begin{array}{l}
\displaystyle{
 v_1(x,y,t) = \frac{i}{k} \Big[ A_f k \cosh(k y) + B_f k \cosh(k(y+H)) + C_f\alpha\cosh(\alpha y) + D_f\alpha\cosh(\alpha(y+H)) \Big] \,\expkw, 
}\medskip\\
\displaystyle{
 v_2(x,y,t) = \Big[ A_f\sinh(k y) + B_f\sinh(k(y+H)) + C_f\sinh(\alpha y) + D_f\sinh(\alpha(y+H)) \Big] \,\expkw , 
}\medskip\\
\displaystyle{
   p(x,y,t) =  \frac{i\rho \omega}{k} \Big[ A_f\cosh(k y) + B_f \cosh(k (y+H)) \Big] \,\expkw ,
}
\end{array}
\label{eq:fluid}
\end{equation} 
where $A_f$, $B_f$, $C_f$ and $D_f$ are constants, and $\alpha$ is given by
\[
\alpha^2=k^2-\frac{i\rho\omega}{\mu}.
\]
Solutions for the displacement of the shell have the form
\begin{equation}
\begin{array}{l}
\displaystyle{
 \us_1(x,t) =-\frac{\mu\theta}{G}  \Big(\frac{i}{k} (B_f k^2 \Sk+D_f \alpha^2 \Sa)+i k (B_f \Sk+D_f \Sa)\Big) \,\expkw ,
} \medskip\\
\displaystyle{
 \us_2(x,t) = \frac{1}{G} \Big( \frac{i\rho\omega}{k} (A_f+B_f \Ck)-2 \mu (A_f k+B_f k \Ck+ C_f \alpha + D_f \alpha  \Ca)\Big) \,\expkw,
}
\end{array}
\label{eq:solid}
\end{equation}
where
\[
   \Ck=\cosh(kH), \qquad  \Sk=\sinh(k H), \qquad \Ca =\cosh(\alpha H), \qquad \Sa =\sinh(\alpha H),
\]
and
\[
G=\Ks+\Ts k^2 - \rhos\hs\omega^2.
\]
The parameter $\theta$ in the solution for $\us_1$ is taken to be 0 for MP-V1 and 1 for MP-V2. 
The four constants $A_f$, $B_f$, $C_f$ and $D_f$ are determined by the boundary conditions at $y=-H$ and the matching conditions at $y=0$, and satisfy
\begin{equation}
\begin{bmatrix}
 -\Sk & 0 & -\Sa & 0 \\
 k \Ck &  k & \alpha \Ca & \alpha \\
\xi & -\Sk G k+ \xi\Ck & 2 i \omega \mu k \alpha & -\Sa G k+2 i \omega \mu k \alpha \Ca \\
  k & k \Ck-2 i \omega \mu\theta k^2 \Sk/G &  \alpha & \alpha \Ca-i \omega \mu\theta(\alpha^2+ k^2) \Sa/G 
\end{bmatrix}
\begin{bmatrix} A_f \\ B_f \\ C_f \\ D_f \end{bmatrix}
= 0 , 
\label{eq:system}
\end{equation}
where $\xi=\rho\omega^2+2 i \omega \mu k^2$. Nontrivial solutions exist if $\omega$ and $k$ satisfy the dispersion relation
\begin{align}
    F(\omega,k) = \det(M)=0,   \label{eq:dispersionV}
\end{align}
where $M=[m_{ij}]\in\Complex^{4\times4}$ is the coefficient matrix in~\eqref{eq:system}.
Given values for $\omega$ and $k$ satisfying~\eqref{eq:dispersionV}, the constants $A_f$, $B_f$ and $C_f$ can be determined in terms of $D_f$ from
\begin{align*}
 \begin{bmatrix}
    m_{11} & m_{12} & m_{13} \\
    m_{21} & m_{22} & m_{23} \\
    m_{31} & m_{32} & m_{33} 
 \end{bmatrix}
  \begin{bmatrix}  A_f \\ B_f \\ C_f \end{bmatrix}
  = 
  - D_f \begin{bmatrix} m_{14} \\ m_{24} \\ m_{34} \end{bmatrix} .
\end{align*}
The parameter $D_f$ is chosen so that the maximum displacement of the shell is $\ampe$, 
i.e.~$\vert\usv(0,0)\vert=\sqrt{\us_1^2(0,0)+\us_2^2(0,0)}=\ampe$. 
The real and imaginary parts of~\eqref{eq:fluid} and~\eqref{eq:solid} define real solutions to the model problems MP-V1 ($\theta=0$) and  MP-V2 ($\theta=1$).
Values for $\omega$ can be found for given values of $k$, and the other parameters of the problem, using standard numerical root-finding software applied to~\eqref{eq:dispersionV}. 
Some computed values for $\omega$ are provided in Figure~\ref{tab:travelingWaveOmega}.

\bibliographystyle{elsart-num}
\bibliography{journal-ISI,jwb,henshawPapers,henshaw,fsi}

\end{document}